\documentclass[journal,twoside]{IEEEtran}

\usepackage{amsmath,amssymb}
\usepackage{graphicx}

\def\e{\mathop{\rm e}\nolimits}
\def\sgn{\mathop{\rm sgn}\nolimits}
\def\Re{\mathop{\rm Re}\nolimits}
\def\Im{\mathop{\rm Im}\nolimits}

\def\Xint#1{\mathchoice
   {\XXint\displaystyle\textstyle{#1}}%
   {\XXint\textstyle\scriptstyle{#1}}%
   {\XXint\scriptstyle\scriptscriptstyle{#1}}%
   {\XXint\scriptscriptstyle\scriptscriptstyle{#1}}%
   \!\int}
\def\XXint#1#2#3{{\setbox0=\hbox{$#1{#2#3}{\int}$}
     \vcenter{\hbox{$#2#3$}}\kern-.5\wd0}}

\def\dashint{\Xint-}

\font\bb=msbm10 scaled \magstep1 
\def\R{\hbox{\bb R}} 
 
\def\Z{\hbox{\bb Z}}

\numberwithin{equation}{section}

\newtheorem{theorem}{Theorem}{}

\begin{document}


\title{Spectrally Accurate Causality Enforcement using SVD-based Fourier Continuations for High Speed Digital Interconnects}


\bibliographystyle{IEEEtran}

%
%
%

\author{Lyudmyla~L.~Barannyk,~\IEEEmembership{Member,~IEEE,}
Hazem~A.~Aboutaleb,
        Aicha~Elshabini,~\IEEEmembership{Fellow, IEEE \& Fellow, IMAPS,} 
        and~Fred~D.~Barlow,~\IEEEmembership{Senior Member, IEEE   \& Fellow, IMAPS}
\thanks{Lyudmyla~L.~Barannyk is with the Department of Mathematics, University of Idaho, Moscow, ID 83844 USA (e-mail: barannyk@uidaho.edu).}
\thanks{Hazem~A.~Aboutaleb is with the Department of Electrical \& Computer Engineering, University of Idaho, Moscow, ID 83844 USA (e-mail: hazemhassan2@gmail.com).}
\thanks{Aicha~Elshabini and Fred Barlow are with the Department of Electrical \& Computer Engineering, University of Idaho, Moscow, ID 83844 USA (e-mail: elshabini@uidaho.edu; fbarlow@uidaho.edu).}
}


\markboth
{Lyudmyla~L.~Barannyk \MakeLowercase{\textit{et al.}}: Spectrally Accurate Causality Enforcement using SVD-based Fourier Continuations}{IEEE Transactions on Components, Packaging and Manufacturing Technology}%


%

\maketitle


\begin{abstract}
We introduce an accurate and robust technique for accessing causality of network transfer functions given in the form of bandlimited discrete frequency responses.
These transfer functions are commonly used to represent  the electrical response of high speed digital interconnects used on chip and in electronic package assemblies. In some cases small errors in the model development lead to non-causal behavior that does not accurately represent the electrical response and may lead to  a lack of convergence in simulations that utilize these models. 
The approach is based on Hilbert transform relations or Kramers-Kr\"onig dispersion relations and  a construction  of causal Fourier continuations using a regularized singular value decomposition (SVD) method. Given a transfer  function, non-periodic in general, this procedure constructs highly accurate Fourier series approximations on the given frequency interval by allowing the function to be periodic in an extended domain. The causality dispersion relations are enforced spectrally and exactly. This eliminates the necessity of approximating the transfer function behavior at infinity and explicit computation of the Hilbert transform. 
We perform the error analysis of the method and take into account a possible presence of a noise or approximation errors in data. 
The developed error estimates can be used in verifying causality of the given data.
The performance of the method is tested on several analytic and simulated examples that demonstrate  an excellent accuracy and reliability of the proposed technique in agreement with the obtained error estimates. 
The method is capable of detecting very small localized causality violations with amplitudes close to the machine precision.

\end{abstract}


\IEEEpeerreviewmaketitle

\section{Introduction} \label{Introduction}

The design of high speed interconnects, that are common on chip and at the package level in digital systems, requires systematic simulations at different levels in order to evaluate the overall electrical system performance and avoid signal and power integrity problems \cite{Swaminathan_Engin_2007}. To conduct such simulations, one needs suitable models that capture the relevant electromagnetic phenomena that affect the signal and power quality. These models are often obtained either from direct measurements or electromagnetic simulations in the form of discrete port frequency responses that represent scattering, impedance, or admittance transfer functions or transfer matrices in scalar or multidimensional cases, respectively. Once frequency responses are available, a corresponding macromodel can be derived using several techniques such as the Vector Fitting  \cite{Gustavsen_Semlyen_1999},  the Orthonormal Vector Fitting \cite{Deschrijver_Haegeman_Dhaene_2007}, the Delay Extraction-Based Passive Macromodeling \cite{Charest_Achar_Nakhla_Erdin_2009} among others. However, if the data are contaminated by errors, it may not be possible to derive a good model. These errors may be due to a noise, inadequate calibration techniques or imperfections of the test set-up in case of direct measurements or approximation errors due to the meshing techniques, discretization errors and errors due to finite precision arithmetic occurring in numerical simulations. Besides, these data are typically available over a finite frequency range as discrete sets with a limited number of samples. All this may affect the performance of the macromodeling algorithm resulting in non-convergence or inaccurate models. Often the underlying cause of such behavior is the lack of causality in a given set of frequency responses \cite{Triverio_Grivet_Talocia_Nakhla_Canavero_Achar_2007}.

Causality can be characterized either in the time domain or the frequency domain. In the time domain, a system is said to be causal if the effect always follows the cause. This implies that a time domain impulse response function $h(t)=0$ for $t<0$, 
and a causality violation is stated if any nonzero value of $h(t)$
is found 
for some $t<0$. 
To analyze causality, one can convert the  frequency responses to the time domain
using the inverse discrete Fourier transform. This approach suffers from the well known Gibbs phenomenon that is inherent for functions that are not smooth enough and represented by a truncated Fourier series. Examples of such functions include  impulse response functions of typical interconnects that have jump discontinuities and whose spectrum is truncated since the  frequency response data are available only on a finite length frequency interval. Direct application of the inverse discrete Fourier transform to raw frequency response data causes severe over and under shooting near the singularities. 
This problem is usually addressed by windowing the Fourier data to deal with the slow decay of the Fourier spectrum \cite[Ch. 7]{Oppenheim_Schafer_1989}. Windowing can also be applied in the Laplace domain \cite{Blais_Cimmino_Ross_Granger_2009} to respect causality. This approach is shown to be more accurate and efficient than the Fourier approach  \cite{Granger_Ross_2009}.
There are other 
filtering techniques that deal with the Gibbs phenomenon but they require some knowledge of location of singularities (see \cite{Gottlieb_Shu_1997, Gelb_Tanner_2006, Tadmor_2007, Mhaskar_Prestin_2009} and references therein). A related paper \cite{Beylkin_Monzon_2009} employs nonlinear extrapolation of Fourier data to avoid the Gibbs phenomenon and the use of windows/filtering.

In the frequency domain, a system is said to be causal if a frequency response given by the transfer function $H(w)$ satisfies the dispersion relations also known as Kramers-Kr\"onig relations \cite{Kramers_1927, Kronig_1926}.  
The dispersion relations can be written using the Hilbert transform. 
They represent the fact that the real and imaginary parts of a causal function are related through Hilbert transform. The Hilbert transform may be expressed in both continuous and discrete forms and is widely used in circuit analysis, digital signal processing, remote sensing and image reconstruction \cite{Guillemin_1977, Oppenheim_Schafer_1989}. 
Applications in electronics include reconstruction \cite{Amari_Gimersky_Bornemann_1995} and correction \cite{Tesche_1992} of measured data, delay extraction \cite{Knockaert_Dhaene_2008}, interpolation/extrapolation of frequency responses \cite{Narayana_Rao_Adve_Sarker_Vannicola_Wicks_Scott_1996}, 
time-domain conversion \cite{Luo_Chen_2005}, estimation of optimal bandwidth and data density using causality checking \cite{Young_2010} and causality enforcement techniques using generalized dispersion relations \cite{Triverio_Grivet_Talocia_2006, Triverio_Grivet_Talocia_2006_2, Triverio_Grivet_Talocia_2008}, causality enforcement using minimum phase and all-pass decomposition and delay extraction \cite{Mandrekar_Swaminathan_2005_2, Mandrekar_Swaminathan_2005, Mandrekar_Srinivasan_Engin_Swaminathan_2006, Lalgudi_Srinivasan_Casinovi_Mandrekar_Engin_Swaminathan_Kretchmer_2006},  causality verification using  minimum phase and all-pass decomposition  that avoids Gibbs errors \cite{Xu_Zeng_He_Han_2006}, causality characterization through analytic continuation for $L_2$ integrable functions \cite{Dienstfrey_Greengard_2001}, causality enforcement using periodic polynomial continuations \cite{Aboutaleb_Barannyk_Elshabini_Barlow_WMED13, Barannyk_Aboutaleb_Elshabini_Barlow_IMAPS, Barannyk_Aboutaleb_Elshabini_Barlow_IMAPS2014} and the subject of the current paper.

The Hilbert transform that relates the real and imaginary parts of a transfer function $H(w)$ is defined on the infinite domain which can be reduced to $[0,\infty)$ by symmetry properties of $H(w)$ for real impulse response functions. However, the frequency responses are usually available over a finite length frequency interval, so the infinite domain is either truncated or behavior of the function for large $w$ is approximated. This is necessary since measurements can only be practically conducted over a finite frequency range and often the cost of the measurements scales in an exponential manner with respect to frequency. Likewise simulation tools have a limited bandwidth and there is a computational cost associated with each frequency data point that generally precludes very large bandwidths in these data sets. Usually $H(w)$ is assumed to be square integrable, which would require the function to decay at infinity. When a function does not decay at infinity or even grows, generalized dispersion relations with subtractions can be successfully used to reduce the dependence on high frequencies and  allow a domain truncation \cite{Triverio_Grivet_Talocia_2006, Triverio_Grivet_Talocia_2006_2, Triverio_Grivet_Talocia_2008}. A review of some previous work on generalized dispersion relations and other methods that address the problem of having finite frequency range is provided in \cite{Triverio_Grivet_Talocia_2008}. 

We take another approach and instead of approximating the behavior of $H(w)$ for large $w$ or truncating the domain, we construct a causal periodic continuation  or causal Fourier continuation of $H(w)$ by requiring the transfer function to be periodic and causal in an extended domain of finite length.
In  \cite{Aboutaleb_Barannyk_Elshabini_Barlow_WMED13, Barannyk_Aboutaleb_Elshabini_Barlow_IMAPS, Barannyk_Aboutaleb_Elshabini_Barlow_IMAPS2014}, polynomial periodic continuations were used to make a transfer function periodic on an extended frequency interval. In these papers, the raw frequency responses were used on the original frequency interval. Once a periodic continuation is constructed, the spectrally accurate Fast Fourier Transform \cite{Cooley_Tukey_1965} implemented in FFT/IFFT routines can be used to compute discrete Hilbert transform and enforce causality. The accuracy of the method was shown to depend primarily on the degree of the polynomial, which implied the smoothness  up to some order of the continuation at the end points of the given frequency domain. This in turn allowed to reduce the boundary artifacts compared to applying the discrete Hilbert transform directly to the data without any periodic continuation, which is implemented in  the function {\tt hilbert} from the popular software Matlab.  
 
In the current work we implement the idea of periodic continuations of an interconnect transfer function by approximating this function  with a causal Fourier series in an extended domain. The approach allows one to obtain extremely accurate approximations of the given function on the original interval. The causality conditions are imposed exactly and directly on Fourier coefficients, so there is no need to compute Hilbert transform  numerically. This eliminates the necessity of approximating the behavior of the transfer function at infinity similar to polynomial continuations employed in \cite{Aboutaleb_Barannyk_Elshabini_Barlow_WMED13, Barannyk_Aboutaleb_Elshabini_Barlow_IMAPS, Barannyk_Aboutaleb_Elshabini_Barlow_IMAPS2014}, and does not require the use of Fast Fourier Transform.  
The advantage of the method is that it is capable of detecting very small localized causality violations with amplitude close to the machine precision, at the order of $10^{-13}$, and a small uniform approximation error can be achieved on the entire original frequency interval, 
so it does not have boundary artifacts reported by using {\tt hilbert} Matlab function,  polynomial continuations \cite{Aboutaleb_Barannyk_Elshabini_Barlow_WMED13, Barannyk_Aboutaleb_Elshabini_Barlow_IMAPS, Barannyk_Aboutaleb_Elshabini_Barlow_IMAPS2014} or generalized dispersion relations \cite{Triverio_Grivet_Talocia_2006, Triverio_Grivet_Talocia_2006_2, Triverio_Grivet_Talocia_2008}. 
%
The performed error analysis unbiases an error due to approximation of a transfer function with a causal Fourier series from causality violations that are due to the presence of a noise or approximation errors in data. 
%
%
The developed estimates of upper bounds for these errors  can be used in checking causality of the given data.
%

The paper is organized as follows. Section \ref{causality_dispersion_relations} provides a background on causality for linear time-translation invariant systems, dispersion relations and  the motivation for the proposed method. In Section \ref{Fourier_continuation} we derive causal spectrally accurate Fourier continuations using truncated singular value decomposition (SVD). 
In Section \ref{error_analysis} we perform the error analysis of the method and take into account a possible noise or approximation errors in the given data. We outline an approach for verifying causality of the given data by using the developed  error estimates. 
In Section \ref{examples}, the technique is applied to several analytic and simulated examples, both causal and non-causal, to show the excellent performance of the proposed method that works in a very good agreement with the developed error estimates. 
Finally, in Section \ref{conclusions} we present our conclusions.

\section{Causality for Linear Time-Translation Invariant Systems} \label{causality_dispersion_relations}

Consider a linear and time-invariant physical system with the impulse response  ${\mathbf h}(t,t')$ subject to a time-dependent input ${\mathbf f(t)}$, to which it responds by an output ${\mathbf x(t)}$.  Linearity of the system implies that the output ${\mathbf x(t)}$ is a linear functional of the input ${\mathbf f(t)}$, while 
time-translation invariance means that if the input is shifted by some time interval $\tau$, the output is also shifted by the same interval, and, hence, the impulse response function ${\mathbf h(t,t')}$  depends only on the difference between the arguments. Thus, the response ${\mathbf x(t)}$ can be written as the convolution of the input ${\mathbf f(t)}$ and the impulse response  ${\mathbf h(t-t')}$ \cite{Nussenzveig_1972}
\begin{equation} \label{1_3_2}
{\mathbf x(t)}=\int_{-\infty}^\infty {\mathbf h(t-t')} {\mathbf f(t')}dt' = {\mathbf h(t)}* {\mathbf f}(t).
\end{equation}
Denote by 
\begin{equation} \label{1_3_3}
{\mathbf H}(w)=\int_{-\infty}^\infty {\mathbf h}(\tau)\e^{-i w\tau}d\tau
\end{equation}
the Fourier transform of ${\mathbf h(t)}$\footnote{Please note that we use an opposite sign of the exponent in the definition of the Fourier transform than in \cite{Nussenzveig_1972}.}. ${\mathbf H}(w)$ is also called the transfer matrix in multidimensional case or transfer function  in a scalar case.

The system is causal if  the output cannot precede the input, i.e. if ${\mathbf f
}(t)=0$  for $t<T$, the same must be true for ${\mathbf x(t)}$. This primitive causality condition in the time domain  implies
%
${\mathbf h}(\tau)=0$,  $\tau<0$,
%
and (\ref{1_3_3}) becomes
\begin{equation} \label{1_3_6}
{\mathbf H}(w)=\int_{0}^\infty {\mathbf h}(\tau)\e^{-i w\tau}d\tau.
\end{equation}
Note that the integral in (\ref{1_3_6}) is extended only over a half-axis, which implies that ${\mathbf H}(w)$ has a regular analytic continuation in lower half $w$-plane.

Examples of physical systems that satisfy the above conditions include electric networks with ${\mathbf f}$ the input voltage, ${\mathbf x}$ the output current,  ${\mathbf H}(w)$ the admittance of the network; ${\mathbf f}$ the input current, ${\mathbf x}$ the output voltage,  ${\mathbf H}(w)$ the impedance; ${\mathbf f}$, ${\mathbf x}$ both  power waves, ${\mathbf H}(w)$ the scattering.

For simplicity, we consider the case with a scalar impulse response $h(t)$ but the approach can also be extended to the multidimensional case for any element of the impulse response matrix ${\mathbf h(t)}$. 

Very often it is assumed that $H(w)$ is square integrable \cite{Dienstfrey_Greengard_2001, Knockaert_Dhaene_2008}, i.e.
%
$\int_{0}^\infty |H(w)|^2 dw < C$
%
for some constant $C$. Then one can use Parseval's theorem to show that  $h(t)$ is also square integrable \cite{Nussenzveig_1972}. The converse also holds \cite{Dym_McKean_1985}. Square integrability of $H(w)$ is often related with the requirement that the total energy of the system is  finite.  Starting from Cauchy's theorem and using contour integration, one can show \cite{Nussenzveig_1972} that for any point $w$ on the real axis, $H(w)$ can be written as
\begin{equation} \label{1_6_7}
H(w)=\frac{1}{\pi i}\dashint_{-\infty}^\infty \frac{H(w')}{w-w'}dw', \quad \mbox{real} \ w,
\end{equation}
where 
\[
\dashint_{-\infty}^\infty=P\int_{-\infty}^\infty = \lim_{\epsilon\to 0}\left(\int_{-\infty}^{w-\epsilon}+\int_{w+\epsilon}^{\infty}\right)
\]
denotes Cauchy's principal value. Separating the real and imaginary parts of (\ref{1_6_7}), we get
\begin{equation} \label{1_6_10}
\Re H(w)=\frac{1}{\pi}\dashint_{-\infty}^\infty \frac{\Im H(w')}{w-w'}dw',
\end{equation}
\begin{equation} \label{1_6_11}
\Im H(w)=-\frac{1}{\pi}\dashint_{-\infty}^\infty \frac{\Re H(w')}{w-w'}dw'.
\end{equation}
These expressions relating  $\Re H$ and $\Im H$ are called the dispersion relations or Kramers-Kr\"onig  relations after Kr\"onig \cite{Kronig_1926} and Kramers \cite{Kramers_1927} who derived the first known dispersion relation for a causal system of a dispersive medium. In mathematics the dispersion relations (\ref{1_6_10}), (\ref{1_6_11}) are also known as the  Sokhotski--Plemelj  formulas. These formulas show that $\Re H$ at one frequency is related to $\Im H$ for all frequencies, and vice versa. Choosing either $\Re H$ or $\Im H$ as an arbitrary square integrable function, then the other one is   completely determined by causality. Recalling that the Hilbert transform is defined
\[
{\mathcal H}[u(w)]=\frac{1}{\pi}\dashint_{-\infty}^\infty \frac{u(w')}{w-w'} dw',
\]
we see that $\Re H$ and $\Im H$ are Hilbert transforms of each other, i.e.
\[
\Re H(w)={\mathcal H}[\Im H(w)], \quad \Im H(w)=-{\mathcal H}[\Re H(w)].
\]
For example, a function $H(w)=\frac{1}{w-i}$ is clearly square integrable and satisfies the dispersion relations (\ref{1_6_10}),  (\ref{1_6_11}), which can be verified by contour integration. 
An example of a function  $H(w)$ that is not square integrable but satisfies the Kramers-Kr\"onig dispersion relations  (\ref{1_6_10}),  (\ref{1_6_11}) is provided by $H(w)=\e^{-iaw}$, $a>0$. The real and imaginary parts are $\cos(aw)$ and $-\sin(aw)$, and dispersion relations  (\ref{1_6_10}),  (\ref{1_6_11}) can be verified by noting that ${\mathcal H}[\cos(aw)]=\sin(aw)$ and ${\mathcal H}[\sin(aw)]=-\cos(aw)$.
%

In practice, the function $H(w)$ may not satisfy the assumption of square integrability and it may be only  bounded or even behave like $O(w^n)$, when $|w|\to\infty$, $n=0,1,2,\ldots$. In such cases, instead of dispersion relations (\ref{1_6_10}), (\ref{1_6_11}), one can use generalized dispersion relations with subtractions, in which a square integrable  function is constructed by  subtracting a Taylor polynomial of $H(w)$ around $w=w_0$ from $H(w)$ and dividing the result by $(w-w_0)^n$. This approach makes the integrand in the generalized dispersion relations less dependent on the high-frequency behavior of $H(w)$. This may be very beneficial when the high-frequency behavior of $H(w)$ may not be known with sufficient accuracy or not accessible in practice at all due to availability of only finite bandwidth data. The technique was proposed in \cite{Beltrami_Wohlers_1966, Nussenzveig_1972} and implemented successfully in \cite{Triverio_Grivet_Talocia_2006, Triverio_Grivet_Talocia_2006_2, Triverio_Grivet_Talocia_2008} to reduce sensitivity of Kramers-Kr\"onig dispersion relations  (\ref{1_6_10}),  (\ref{1_6_11}) to the high-frequency data. 

In this paper, we take an alternative approach motivated by the example of the periodic function $H(w)=\e^{-iaw}$, $a>0$, mentioned above, that is not square integrable but still satisfies Kramers-Kr\"onig dispersion relations  (\ref{1_6_10}),  (\ref{1_6_11}). 
The transfer function $H(w)$ in practice is typically known only over a finite frequency interval with the limited number of discrete values and it is not periodic in general.
Direct application of  dispersion relations  (\ref{1_6_10}),  (\ref{1_6_11}) produces large errors in the boundary regions mainly because the high-frequency behavior of $H(w)$ is missing, unless data decay to zero at the boundary. To overcome this problem, we construct a spectrally accurate causal periodic continuation of $H(w)$ in an extended domain. A method for constructing a periodic continuation, also known as Fourier continuation or Fourier extension, which is based on  regularized singular value decomposition (SVD), was recently proposed in  \cite{Boyd_2002, Bruno_2003, Bruno_Han_Pohlman_2007, Lyon_2012} (see also references therein). This method allows one to calculate Fourier series approximations of non-periodic functions such that a Fourier series is periodic in an extended domain. 
Causality can be imposed directly on the Fourier coefficients producing a causal Fourier continuation, thus satisfying causality exactly. 
%
%
The Fourier coefficients  are determined by solving an overdetermined and regularized least squares problem since the system suffers from numerical ill-conditioning. 
%
The resulting causal Fourier continuation is then compared with the given discrete data on the original bandwidth of interest. A decision about causality of the given data is made using the error estimates developed in Section \ref{error_analysis}. 
%

In the next section we provide details of the derivation of causal Fourier continuations.


\section{Causal Fourier Continuation} \label{Fourier_continuation}

Consider a transfer function $H(w)$ available at a set of discrete frequencies from $[w_{min},w_{max}]$, where $w_{min}\geq 0$. First, let $w_{min}=0$, so we have the baseband case. Since equations (\ref{1_6_10}),  (\ref{1_6_11})  are homogeneous in the frequency variable, we can rescale $[0,w_{max}]$ to  $[0,0.5]$ using the transformation $x=\frac{0.5}{w_{max}}w$ 
for convenience, to get a rescaled transfer function $H(x)$.   
The time domain impulse response function $h(t)$ is often real-valued. Hence, the real and imaginary parts of  $H(w)$, as the Fourier transform of $h(t)$,  and, hence, of $H(x)$, are even and odd functions. This implies that the  discrete set of  rescaled frequency responses $H(x)$  is available on the unit length  interval  $x\in[-0.5,0.5]$ by spectrum symmetry. In some cases, the data are available only from a non-zero, low-frequency cutoff  $w_{min}>0$, which corresponds to the bandpass case. 
The proposed procedure is still applicable since it does not require data points to be equally spaced. The transmission line example \ref{transmission_line_example} considers such situation.

The idea is to construct an accurate Fourier series approximation of $H(x)$ by allowing the Fourier series to be periodic and causal in an extended domain. The result is the Fourier continuation of $H$ that we denote by ${\mathcal C}(H)$, and it is defined by
\begin{equation}\label{E1}
{\mathcal C}(H)(x)=\sum_{k=-M+1}^{M} \alpha_k \e^{-\frac{2\pi i}{b} k x},
\end{equation}
for even number $2M$ of terms, whereas for odd number $2M+1$ of terms, the index $k$ varies from $-M$ to $M$. 
Throughout this paper we will consider Fourier series with even number of terms for simplicity. All presented results have analogues for Fourier series with odd number of terms. Here $b$ is the period of approximation. For SVD-based periodic continuations $b$ is normally chosen as twice the length of the domain on which function $H$ is given \cite{Bruno_Han_Pohlman_2007}. The value $b=2$ is not necessarily optimal and it is shown \cite{Lyon_2012a} to depend  on a function being approximated. 
%
%
For causal Fourier continuations we also find that the optimal value of $b$ depends on $H$. In practice, $b$ can be varied in $1<b\leq 4$  to get more optimal performance of the method. For very smooth functions, it is better to use a wider extension zone with $b\geq 2$, for example, $b=2$ or $b=4$ was enough in most of our examples. However, for functions that are wildly oscillatory or have high gradients in the boundary regions of the domain where the original data are available, a smaller extension zone with $1<b<2$ is recommended   \cite{Boyd_2002}.   We used  $b=1.1$ in one of our examples. 
%
Assume that values of the function $H(x)$ are known at $N$ discretization or collocation points $\{x_j\}$, $j=1,\ldots,N$, $x_j\in[-0.5,0.5]$.  Note that  ${\mathcal C}(H)(x)$ is a trigonometric polynomial of degree at most $M$.
 
Since $\Re H(x)$ and $\Im H(x)$ are even and odd functions of $x$, respectively, the Fourier coefficients 
\[
\alpha_k=\frac{1}{b}\int_{-b/2}^{b/2} H(x) \overline{\phi_k(x)}dx, \quad k=1,\ldots, M,
\]
are real. Here $\phi_{k}(x)=\e^{-\frac{2\pi i }{b}k x}$, $k\in\Z$, and $\bar{ \  }$ denotes the complex conjugate. Functions  $\{\phi_{k}(x)\}$ form a complete orthogonal basis in $L_2[-\frac b2, \frac b2]$, and, in particular
\begin{equation} \label{Kronecker_delta_basis}
\int_{-b/2}^{b/2} \phi_k(x)\overline{\phi_{k'}(x)} dx = b \,\delta_{kk'},
\end{equation}
where $\delta_{kk'}$ is the Kronecker delta. In addition, $\overline{\phi_k}(x)=\e^{\frac{2\pi i }{b}k x}=\phi_{-k}(x)$.

For a function $\e^{-iax}$, the Hilbert transform is
$
{\mathcal H}\{\e^{-iax}\}=  i\sgn(a) \e^{-iax}
$.
%
Hence, 
%
\begin{equation}\label{E2}
{\mathcal H}\{\phi_k(x)\}= i\sgn(k)  \phi_k(x),
\end{equation}
which implies that the  functions $\{\phi_k(x)\}$ are the eigenfunctions of the Hilbert transform ${\mathcal H}$ with associated eigenvalues $\pm i$ with $x\in[-\frac b2,\frac b2]$. We will use relations (\ref{E2}) to impose a causality condition on the coefficients of $C(H)(x)$ similarly as it was done in \cite{Knockaert_Dhaene_2008} for the case of square integrable $H(w)$ where the idea of projecting on the eigenfunctions of the Hilbert transform in $L_2(\R)$ \cite{Weideman_1995} was used. In the present work, the 
square integrability of $H(w)$ is not required and more general transfer functions than in \cite{Knockaert_Dhaene_2008} can be considered.


For convenience of derivation, let us  write $C(H)(x)$ as a Fourier series 
$
{\mathcal C}(H)(x)=\sum_{k=-\infty}^{\infty} \alpha_k \phi_k (x)
$,
which will be truncated  at the end to get a Fourier continuation in the form (\ref{E1}).
Let ${\mathcal C}(H)(x)=\Re{\mathcal C}(H)(x) + i \Im {\mathcal C}(H)(x)$ and $\phi_k(x) =\Re \phi_k(x)+i \Im \phi_k(x)$. 
%
%
%
%
%
%
Since
\[
\Re \phi_k=\frac 12 (\phi_k+\overline{\phi}_k), \quad \Im \phi_k=\frac{1}{2i} (\phi_k-\overline{\phi}_k)
\]
we obtain
\[
\Re{\mathcal C}(H)(x)=\sum_{k=-\infty}^{\infty} \alpha_k \Re \phi_k=\frac 12  \sum_{k=-\infty}^{\infty} \alpha_k  (\phi_k+\overline{\phi}_k)
\]
and, since $\overline{\phi_k}=\phi_{-k}$, we have
\[
\Re{\mathcal C}(H)(x)=\frac 12 \sum_{k=-\infty}^{\infty} \alpha_k   (\phi_k+\phi_{-k}) = 
\frac 12 \sum_{k=-\infty}^{\infty} (\alpha_k +\alpha_{-k})  \phi_k,
\]
where in the last sum we changed the order of summation in the second term.
Similarly, we can show that
%
%
\[
\Im{\mathcal C}(H)(x)=\frac {1}{2i}  \sum_{k=-\infty}^{\infty}  (\alpha_k -\alpha_{-k})  \phi_k.
\]

For a causal periodic continuation, we need $\Im{\mathcal C}(H)(x)$ to be the  Hilbert transform of $-\Re{\mathcal C}(H)(x)$. Hence,
\[
\frac {1}{2i}  \sum_{k=-\infty}^{\infty}  (\alpha_k -\alpha_{-k})  \phi_k= -{\mathcal H}\left[\frac 12 \sum_{k=-\infty}^{\infty} (\alpha_k +\alpha_{-k})  \phi_k \right].
\] 
Employing linearity of the Hilbert transform, we get
\[
\frac {1}{2i}  \sum_{k=-\infty}^{\infty} (\alpha_k -\alpha_{-k})  \phi_k= -\frac 12 \sum_{k=-\infty}^{\infty}  (\alpha_k +\alpha_{-k})  {\mathcal H}[\phi_k].
\]
Using (\ref{E2}), we obtain 
\[
\frac {1}{2i}  (\alpha_k -\alpha_{-k}) =-\frac 12 (\alpha_k +\alpha_{-k})  i \sgn(k) \quad \mbox{for any} \ k\in\Z
\]
or
\[
\alpha_k-\alpha_{-k} = (\alpha_k +\alpha_{-k}) \sgn(k), \quad k\in\Z, 
\]
that implies  $\alpha_{k}=0$ for $k\leq 0$ in  (\ref{E1}).
Hence, a causal Fourier continuation has the form
\begin{equation} \label{E3_0}
{\mathcal C}(H)(x)=\sum_{k=1}^{M}  \alpha_k \phi_k(x) 
\end{equation}
where we truncated the infinite sum to obtain a trigonometric polynomial.
Evaluating $H(x)$ at points $x_j$, $j=1,\ldots,N$, $x_j\in[-0.5, 0.5]$, produces
a complex valued system
\begin{equation} \label{E3}
{\mathcal C}(H)(x_j) =\sum_{k=1}^{M}  \alpha_k \phi_k(x_j)
\end{equation}
with $N$ equations for $M$ unknowns $\alpha_{k}$, $k=1,\ldots,M$, $N\geq M$. If $N>M$,  the system (\ref{E3})  is overdetermined and has to be solved in the least squares sense. 
%
%
When Fourier coefficients $\alpha_k$ are computed, formula (\ref{E3_0}) provides reconstruction of $H(x)$ on $[-0.5, 0.5]$. 


To ensure that  numerically computed Fourier coefficients $\alpha_k$ are real, instead of solving complex-valued system (\ref{E3}), one  can separate the real and imaginary parts of ${\mathcal C}(H)(x_j)$ to obtain real-valued system
\begin{equation} \label{E4}
\begin{array}{l}
\displaystyle
\phantom{-}
\Re{\mathcal C}(H)(x_j)=\sum_{k=1}^{M}  \alpha_{k}   \Re\phi_k(x_j), \\[13pt]
\displaystyle
\hspace{10pt}
\Im{\mathcal C}(H)(x_j)=\sum_{k=1}^{M}  \alpha_{k} \Im\phi_k(x_j).
\end{array}
\end{equation}
%
This produces $2N$ equations ($N$ equations for both real and imaginary parts) and $M$ unknowns $\alpha_{k}$. We show below that both complex  (\ref{E3}) and real  (\ref{E4})  formulations give the reconstruction errors of the same order with the real formulation performing slightly better. To distinguish between the continuation ${\mathcal C}(H)$ computed using complex  or real formulation, we will use notation ${\mathcal C}^C(H)$ and ${\mathcal C}^R(H)$, respectively.

Consider the real formulation  (\ref{E4}) and introduce the following notation. Let $\vec f=\bigl(\Re H(x_1), \ldots, \Re H(x_N), \Im H(x_1), \ldots, \Im H(x_N)\bigr)^T$, 
$\vec\alpha=\bigl(\alpha_{1},\ldots,\alpha_{M}\bigr)^T$, where ${}^T$ denotes the transpose,  and  matrix $A$ have entries
\[
\displaystyle
A_{jk} =\Re\{\e^{-\frac{2\pi i}{b}k x_j}\},
\ j=1,\ldots,N, \ k=1,\ldots, M,
\] \[
A_{(j+N),k}
=\Im\{\e^{-\frac{2\pi i}{b}k x_j}\},
\  j=1,\ldots,N, \ k=1,\ldots, M.
\]
%
Similar notation can be made for the complex formulation (\ref{E3}). 
Then the coefficients $\alpha_{k}$, $k=1,\ldots, M$, are defined as a least squares solution of 
%
$
A\vec\alpha=\vec f
$
%
written as
\[
\min_{\{\alpha_k\}} \sum_{j=1}^{2N} \left| \sum_{k=1}^{M}   \alpha_{k} A_{jk}-f_j\right|^2,
\]
that minimizes the Euclidean norm of the residual. This least squares problem is extremely ill-conditioned, as explained in \cite{Huybrechs_2010} using the theory of frames. However, it can be regularized using a truncated SVD method when singular values below some cutoff tolerance $\xi$  close to the machine precision are being discarded. 
%
%
%
In this work we use $\xi=10^{-13}$ as the threshold to filter the singular values.
%
The ill-conditioning increases as $M$ increases by developing rapid  oscillations in the  extended region. These oscillations are typical for SVD-based Fourier continuations. Once the system reaches a critical size that does not depend on the function being approximated,  the coefficient matrix becomes rank deficient and  the regularization of the SVD is required to treat singular values close to the machine precision. Because of the rank deficiency, the Fourier continuation is not longer unique. Applying the truncated SVD method produces the minimum norm solution $\{\alpha_k\}$, $k=1,\ldots, M$,
for which the corresponding Fourier continuation is oscillatory. The oscillations in the extended region do not significantly affect  the quality of the causal Fourier continuation on the original domain and varying $b$ can minimize their effect  and decrease the overall reconstruction error, 
especially in the boundary domain. 

Another way to make ill-conditioning of matrix problems (\ref{E3}) or (\ref{E4}) better is to use more data (collocation) points $N$ than the Fourier coefficients $M$. This is called ``overcollocation" \cite{Boyd_2002} and makes the problem more overdetermined and helps to increase the accuracy of solutions. It is recommended to use at least twice more collocation points $N$ than the Fourier coefficients $M$, i.e. $N=2M$. The convergence can be checked by keeping the number of Fourier coefficients $M$ fixed and increasing the number of collocation points $N$. The overcollocation also helps with filtering out trigonometric interpolants that have very small errors at collocation points $x_j$ but large oscillations between the collocation points \cite{Boyd_2002}. In all our examples we use at least $N=2M$ as an effective way to obtain an 
accurate and reliable approximation of $H(x)$ over the  interval $[-0.5, 0.5]$.


In the multidimensional case when a transfer matrix ${\mathbf H(w)}$ is given, the above procedure can be extended to all elements of the matrix. Computing SVD is  an expensive numerical procedure for large matrices. The operation count to find a least squares solution  of $Ax=b$ using the SVD method with $A$ being $N\times M$, $N\geq M$, matrix, is of the order $O(NM^2+M^3)$ \cite{LAPACK_Users_Guide}.  The actual CPU times for computing SVD, solving a linear system in the least squares sense and constructing causal Fourier continuations for various values of $M$ used in this work, are shown in Table \ref{T_CPU} in Section \ref{examples}.  
%
%
Computational savings can be achieved by noting that in our problem
the matrix $A$ depends only on the location of frequency points at which transfer matrix is evaluated/available, and continuation parameters $M$ and $b$ but does not depend on actual values of ${\mathbf H(w)}$. Having frequency responses at $N$ points, we can fix $N=2M$,  choose $b=2$ as a default value and compute SVD only once prior to verifying causality. Varying $1<b\leq 4$ or $2M<N$  for each element of ${\mathbf H(w)}$ separately, if needed, would require recomputing SVD. 
%
%



\section{Error Analysis} \label{error_analysis}


In this section, we provide an upper bound for the error in  approximation of a given function $H(x)$ 
by its causal Fourier continuation  ${\mathcal C}(H)(x)$. 
The analysis of convergence of Fourier continuation technique based on the truncated SVD method was done by Lyon \cite{Lyon_2012a} using a split real formulation.  
In this work, we extend his results  to causal Fourier continuations. 
The obtained error estimates can be employed to characterize causality of a given set of data. 
%
%
\subsection{Error estimates}

Denote by $\hat H_{M}$ any function of the form  
\begin{equation}\label{Fourier_series_arbitrary}
\hat H_{M}(x)=\sum_{k=1}^{M}  \hat\alpha_k \phi_k(x)
\end{equation}
where $\phi_k(x)=\e^{-\frac{2\pi i}{b}kx}$, $k=1,\ldots,M$, as before. 
Let $A=U \Sigma V^*$ be the full SVD decomposition \cite{Trefethen_Bau_1997} of the matrix $A$ with entries $A_{kj}=\phi_k(x_j)$, $j=1,\ldots,N$, $k=1,\ldots, M$,
%
%
where $U$ is an $N\times N$  unitary matrix, 
$\Sigma$ is an $N\times M$ diagonal  matrix of singular values $\sigma_j$, $j=1,\ldots,p$, $p=\min({N,M})$,  $V$ is an $M \times M$ unitary matrix with entries $V_{kj}$,
 and $V^*$ denotes the complex conjugate transpose of $V$. 

We can prove the following result.
 

\begin{theorem} 

Consider a rescaled transfer function $H(x)$ defined by symmetry on $\Omega=[-0.5,-a]\cup[a, 0.5]$, where $a=0.5\frac{w_{min}}{w_{max}}$, whose values are available at points $x_j\in \Omega$, $j=1,\ldots,N$. Then the error in approximation of $H(x)$, that is known with some error $\varepsilon$, by 
its causal Fourier continuation ${\mathcal C}(H)(x)$ defined in (\ref{E3_0}) on a wider domain $\Omega^c=[-b/2,b/2]$, $b\geq 1$,
has the  upper bound
\[
||H- {\mathcal C}(H+\varepsilon) ||_{L_2(\Omega)} \leq (1+\Lambda_2 \sqrt{N(M-K)}) 
\] 
\begin{equation}\label{combined_bound_noise}
\times \left(|| H-\hat H_{M} || _{L_\infty(\Omega)} + ||\varepsilon||_{L_\infty(\Omega)} \right)
+\Lambda_1 \sqrt{K/b} ||\hat H_{M}||_{L_\infty(\Omega^c)}
\end{equation}
and holds for all functions of the form (\ref{Fourier_series_arbitrary}). Here
\begin{equation}\label{lambda12}
\Lambda_1=\max_{j:\ \sigma_j<\xi}|| v_j(x) ||_{L_2(\Omega)}, \quad
\Lambda_2=\max_{j:\ \sigma_j>\xi}\frac{||v_j(x) ||_{L_2(\Omega)}}{\sigma_j},
\end{equation}
and  functions $v_j(x)=\sum_{k=1}^{M} V_{kj} \phi_k(x)$ are each an up to $M$ term causal Fourier series with coefficients given by the $j$th column of $V$; $K$ denotes the number of singular values that are discarded, i.e. the number of $j$ for which $\sigma_j<\xi$, where $\xi$ is the cut-off tolerance. 

\end{theorem}

\underline{Proof}.
To obtain the  error bound  (\ref{combined_bound_noise}) we use ideas from \cite{Lyon_2012a} but employ a complex formulation and  impose causality on Fourier coefficients.
%
%
The error bound for $||H- {\mathcal C}(H) ||_{L_2(\Omega)} $ is expressed in terms of the error $||H- \hat H_{M}||_{L_\infty(\Omega)} $ in approximation of a function with a causal Fourier series
 and $||\hat H_{M}||_{L_\infty(\Omega^c)} $ for any given causal Fourier series $\hat H_{M}$.  
This requires finding upper bounds for $||{\mathcal C}(H-\hat H_{M})  ||_{L_2(\Omega)} $ and $||\hat H_{M}- {\mathcal C}(\hat H_{M})||_{L_2(\Omega)}  $, that estimate the error due to truncation of singular values and the effect of Fourier continuation on 
 the error in approximation of a function with a causal Fourier series. If function $H$ is known with some error $\epsilon$, its effect is also included in a straight-forward way.
The bound for $||H- \hat H_{M}||_{L_\infty(\Omega)} $  follows from Jackson Theorems \cite{Cheney_2000} that estimate the error in approximation of a periodic function with its $M$ causal Fourier series $\hat H_{M}$
%
%
as a partial case. Indeed, a causal $M$ mode Fourier series can be considered as an $2M$ mode Fourier series whose coefficients with nonpositive indices are zero. Hence, the error in approximating a $b$-periodic function $H$ with $k$ continuous derivatives with 
a causal $M$ mode Fourier series,  has the following upper bound:
\begin{equation}\label{Fourier_error}
|| H-\hat H_{M} || _{L_\infty(\Omega)}\leq \frac\pi 2 \left(\frac b \pi \right)^k \left(\frac {1}{2M}\right)^k ||H^{(k)}||_{L_\infty(\Omega^c)}.
\end{equation}
%
%
 Left and right singular vectors that form columns of matrices $U$ and $V$ are used in the derivation of the  error estimates (\ref{combined_bound_noise}), (\ref{lambda12}) as alternatives to Fourier basis.
\hfill$\square$ 
 
As can be seen from (\ref{lambda12}),  $\Lambda_1$,  $\Lambda_2$  and $K$ depend only on the continuation parameters $N$, $M$, $b$ and $\xi$ as well as location of discrete points $x_j$, and not on the function $H$.  The behavior of $\Lambda_1 \sqrt{K/b}$ and $\Lambda_2 \sqrt{N(M-K)}$ as functions of $M$ can be investigated using direct numerical calculations. We do this for the case of equally spaced points $x_j$, $j=1,\ldots,N$, which is typical in applications, and use $N=2M$ with $b=2$. Other distributions of points $x_j$, values of $b$ and relations between $M$ and $N$ can be analyzed in a similar manner. For example, results for $b=4$ are very similar to the case with $b=2$.
 We find that while $\Lambda_1$ does not change much with $M$ and remains small, the coefficient $\Lambda_1 \sqrt{K/b}$ stays close to the cut-off value $\xi$ for small $M$ and increases at most to $10\xi$ for large $M$. This behavior does not seem to depend on the cut-off value $\xi$  and the results are similar for $\xi$ varying from $10^{-13}$ to $10^{-7}$. The number $K$ of discarded singular values grows with $M$ (provided that singular values above $\xi$ are computed accurately) since the ill-conditioning of the problem increases with $M$. Values of $\Lambda_2$ remain close to $1$ for values $M$  we considered. At the same time, the coefficient $\Lambda_2 \sqrt{N(M-K)}$ grows approximately as $\sqrt{M}$ as $M$ increases. 

Since the error (\ref{Fourier_error}) in approximation of a function with a causal Fourier series decays as ${\mathcal O}(M^{-k})$ and the  coefficient  $\Lambda_2 \sqrt{N(M-K)}$ grows as ${\mathcal O}(M)$, the  error bound part that is due to a causal Fourier series approximation
\[
\epsilon_F\equiv (1+\Lambda_2 \sqrt{N(M-K)}) || H-\hat H_{M} || _{L_\infty(\Omega)} 
\]
decays at least as fast as ${\mathcal O}(M^{-k+1})$. For comparison, the analogous error bound term for Fourier continuations reported in \cite{Lyon_2012a} is on the order of ${\mathcal O}(M^{-k+1/2})$, i.e. a causal Fourier series converges slightly slower  than  a standard Fourier series.

In practice, the smoothness order $k$ of the transfer function $H(x)$ may not be known. In this case, it can be estimated by noting that the error bound $\epsilon_F$ can be written as
\begin{equation} \label{error_epsilon1}
\epsilon_F\sim \tilde C M^{-k+1}.
\end{equation}
Taking natural logarithm of both sides, we get
\begin{equation}\label{error_ln}
\ln \epsilon_F \sim \ln\tilde C+(-k+1)\ln M,
\end{equation}
i.e. $\ln \epsilon_F$ is approximately a linear function of $\ln M$. The values of $\tilde C$ and $k$ can be estimated as follows. Assume that $H$ is known at $N$ frequency points. Usually the number of frequency responses is fixed and it may not be possible to get data with higher resolution. Assume that the errors due to truncation of singular values (term with $\Lambda_1$) and a noise in data  (term with $||\varepsilon ||_{L_\infty(\Omega)}$)  are small, so that the error due to a causal Fourier series approximation is dominant. Let $E_R(x)$ and $E_I(x)$ be reconstruction errors,
%
%
\begin{equation} \label{error_Re}
E_R(x)=\Re H(x) - \Re{\mathcal C}(H)(x),
\end{equation}
\begin{equation}\label{error_Im}
E_I(x)=\Im H(x) - \Im{\mathcal C}(H)(x)
\end{equation}
%
on the original interval $[-0.5,0.5]$. Compute $E_R(x)$ and $E_I(x)$
with $N$, $N/2$, $N/4$ etc. samples, i.e. with $M$, $M/2$, $M/4$ etc. Fourier coefficients. Solve equation (\ref{error_ln}) in the least squares sense to find approximations of $\ln\tilde C$ and $-k+1$, and, hence, to $\tilde C$ and $k$. Then the error term $\epsilon_F$ can be extrapolated to higher values of $M$ using  (\ref{error_epsilon1}) to see if the causal Fourier series approximation error decreases if the number $M$ of Fourier coefficients increases, i.e. resolution increases.
%


The error bound term 
\begin{equation}\label{error_truncation}
\epsilon_T=\Lambda_1 \sqrt{K/b} ||\hat H_{M}||_{L_\infty(\Omega^c)},
\end{equation}
that is due to the truncation of singular values, is
typically small and close to the cut-off value $\xi$ for small $M< 250$ and at most $10\xi$ for $250\leq M\leq 1500$. 
As can be seen from (\ref{error_truncation}), 
$\epsilon_T$ depends on $b$ and the function $H$ being approximated.
The default value $b=2$ may not provide the smallest  error. 
To find a more optimal value of $b$, a few values in  $1<b\leq 4$ may be tried  to determine which one gives smaller overall reconstruction errors. In case of non-causal functions, varying $b$ does not essentially effect the size of reconstruction errors. 
%


The error $\varepsilon$ in data should be known in practice since the error in measurements  or the accuracy of full wave simulations are typically known. The error bound term due to a noise in data
\[
\epsilon_n=(1+\Lambda_2 \sqrt{N(M-K)}) ||\varepsilon||_{L_\infty(\Omega)} 
\]
consists of the norm of the noise amplified by the coefficient $\Lambda_2 \sqrt{N(M-K)})$ that grows as ${\mathcal O}(M)$  as was shown above. In numerical experiments that we conducted the reconstruction errors due to a noise in data seem to level off to the order of  $\varepsilon$ and are not amplified significantly as the resolution increases. This does not contradict the error estimate (\ref{combined_bound_noise}),  (\ref{lambda12}).  The error bounds are not tight and the actual reconstruction errors may be smaller than the error bounds suggest.

\subsection{Causality characterization}

The  error estimate (\ref{combined_bound_noise}),  (\ref{lambda12}) shows that the reconstruction errors $E_R(x)$ and $E_I(x)$ defined in (\ref{error_Re}),  (\ref{error_Im}) can be dominated by either the error due to approximation of a function with its causal Fourier series, which has the upper bound $\epsilon_F$, or the error $\epsilon$ due to a noise or approximation errors in data, which has the upper bound $\epsilon_n$.  If the only errors in data are  round-off errors, then the reconstruction errors will approach or will be bounded by the error due to truncation of singular values, which has the upper bound $\epsilon_T$.
The noise level $\epsilon$  may be known. In case of experimental data,  $\epsilon$ could be around $10^{-3}$ or $10^{-4}$, for example. Data obtained from finite element simulations may be accurate within $10^{-6}$ or $10^{-7}$, for instance, which would correspond to a single precision accuracy. In some cases, the expected accuracy may be even higher.
In this case if the reconstruction errors are higher than  $\epsilon_n$ (in practice $\epsilon$ can be used), then the reconstruction errors may be dominated by a causal Fourier series approximation error with  the upper bound $\epsilon_F$. Determining the smoothness order of  $H$ using the model (\ref{error_epsilon1}) 
with  $N$, $N/2$, $N/4$ etc. samples can give an answer whether the causal Fourier series approximation error can be made smaller by increasing the resolution of the data. If the model (\ref{error_epsilon1}) extrapolated to higher values of $M$, decays with $M$, this implies that the Fourier series approximation error can be decreased by using more frequency samples, then we can state that the dispersion relations are satisfied within error given by $E_R(x)$ and $E_I(x)$ and causality violations, if present, have the order smaller than or at most  the order of reconstruction errors  $E_R(x)$ and $E_I(x)$. In this case, with fixed resolution, the method will not be able to detect these smaller causality violations. The noise level $\epsilon$ may not be known but it can be determined by comparing reconstruction errors  $E_R(x)$ and $E_I(x)$ as the resolution of data increases, if possible, or decreases, since the reconstruction error due to the noise does not significantly depend on the resolution in practice. This means that as the number $N$ of samples increases, the reconstruction errors level off at the value of the noise $\epsilon$. Equivalently, there may be a situation that using $N$, $N/2$, $N/4$ etc. samples (resolution becomes coarser) the reconstruction error does not increase. Instead, it remains at the same order, which would be a noise level $\epsilon$. In this case, the dispersion relations will be satisfied within the error $\epsilon$, which would be the order of causality violations. Please see Example \ref{example_Micron}  of a finite element model of a package where a related situation is discussed.

%

In the next section we employ the proposed causal Fourier continuation based method to several analytic and simulated examples that are causal/non-causal or have imposed causality violations. When a given transfer function $H(x)$ is causal, we expect that the dispersion relations  (\ref{1_6_10}), (\ref{1_6_11}) are satisfied with the accuracy close to machine precision. This is so called an ideal causality test. When there is a causality violation, the method is expected to point to the location of violation or at least develop reconstruction error close to  the order given by the noise in data as  suggested by the error analysis. We  show that the results are in full agreement with the error estimates developed in Section \ref{error_analysis}.

\section{Numerical Experiments: Causality Verification} \label{examples}


\subsection{Two-pole example} \label{Example1}

\begin{figure}[h] \begin{center}
\includegraphics[width=2.4in,angle=0]{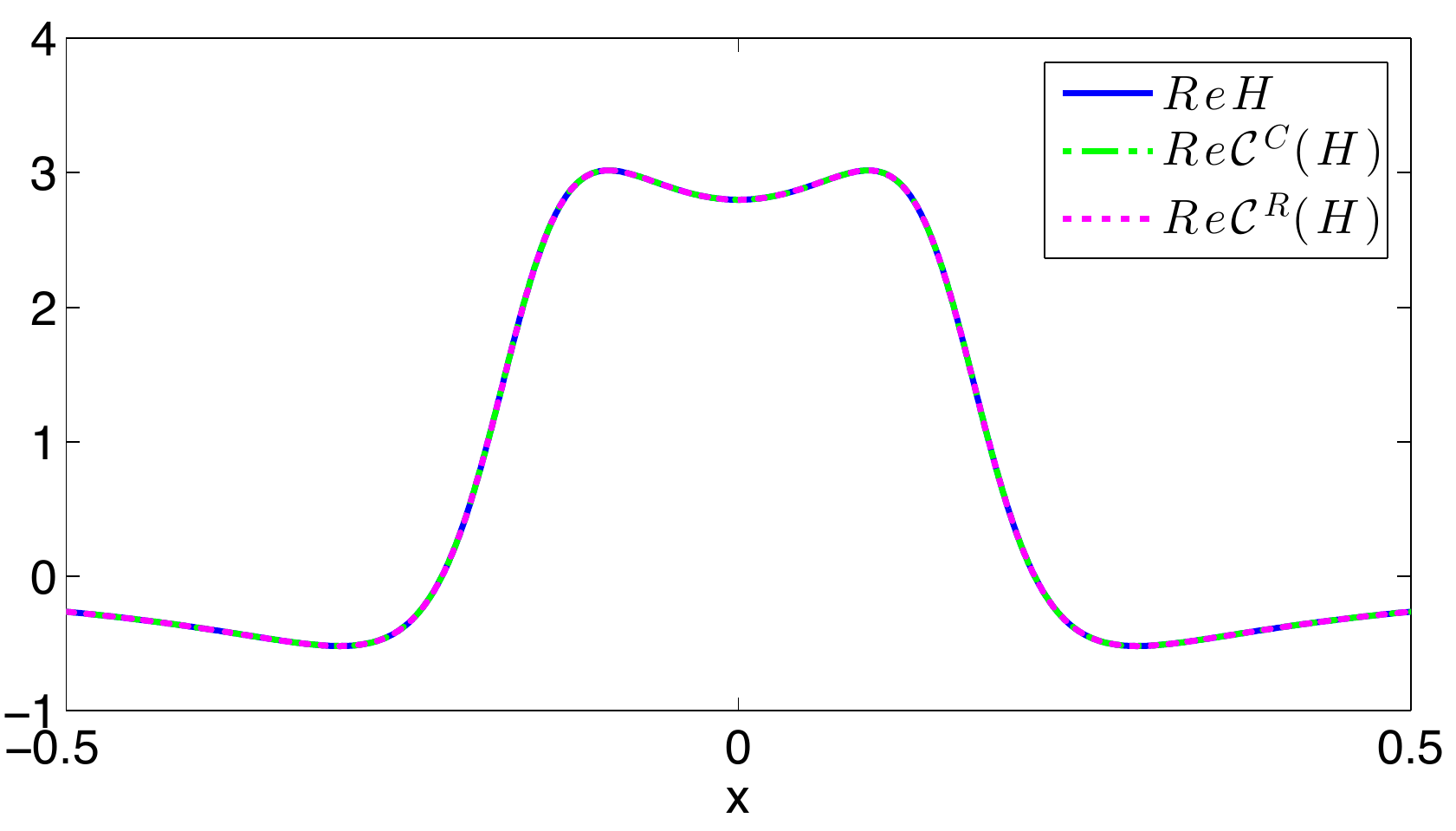}
\includegraphics[width=2.4in,angle=0]{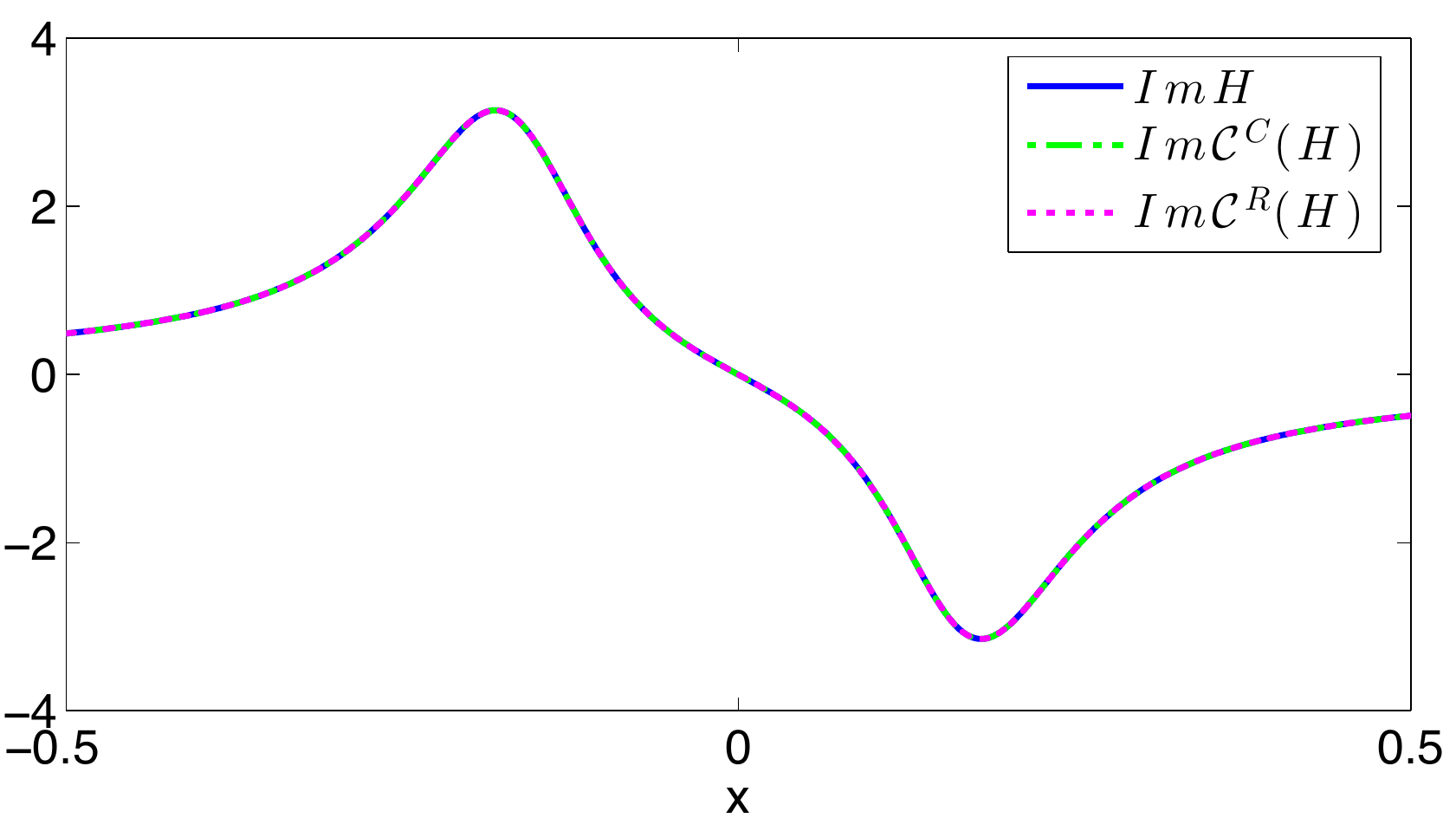}
\end{center}
\caption{$H(x)$ and its Fourier continuations ${\mathcal C}^C(H)$ and ${\mathcal C}^R(H)$ computed using complex  (\ref{E3}) and real  (\ref{E4}) formulations in example \ref{Example1} with $M=250$,  $N=4M$, $b=4$ shown on $[-0.5, 0.5]$.} 
\label{F4}
\end{figure}

\begin{figure}[h] \begin{center}
\includegraphics[width=2.4in,angle=0]{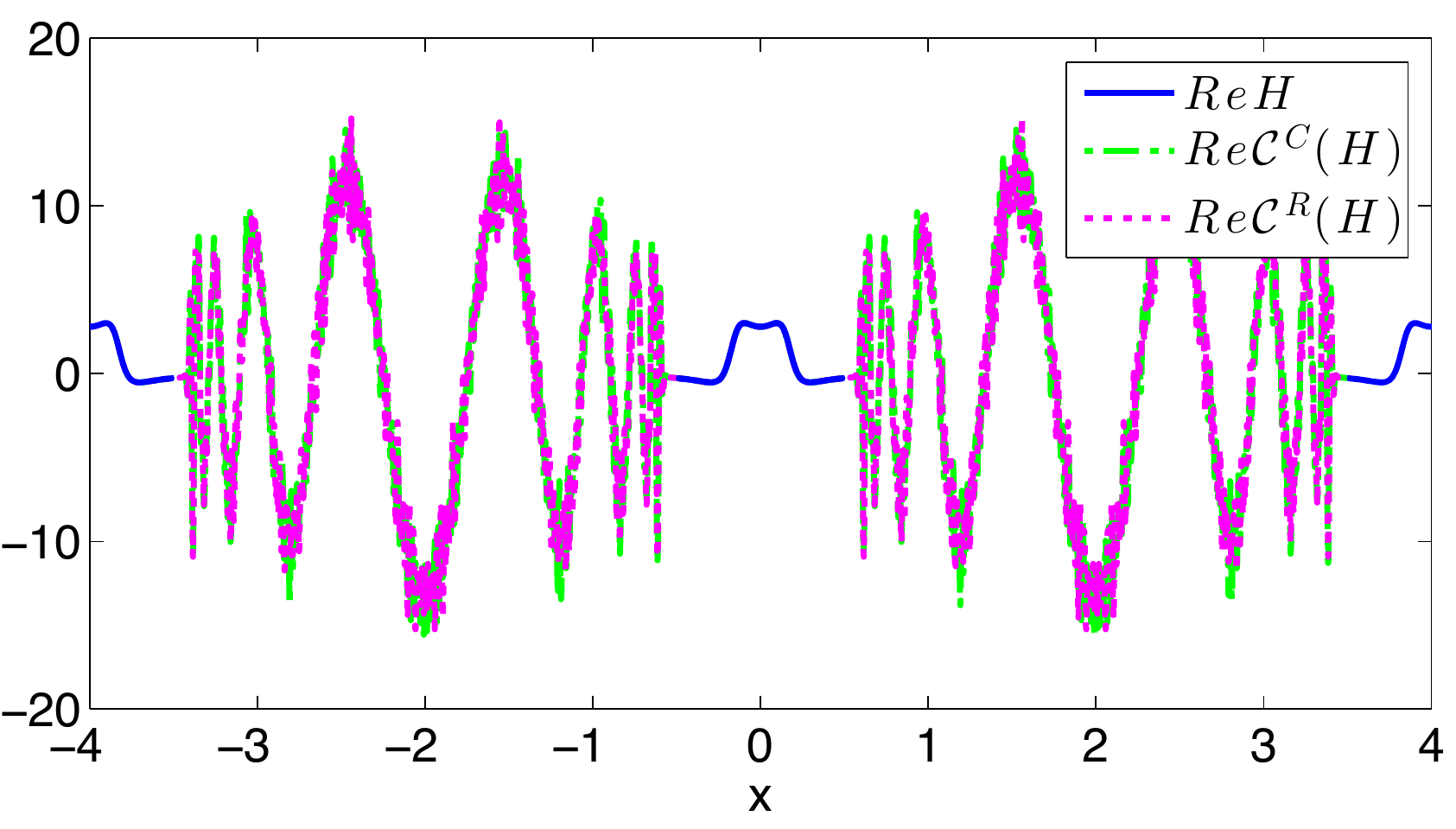}
\end{center}
\caption{The real part of Fourier continuations ${\mathcal C}^C(H)$ and  ${\mathcal C}^R(H)$ in example \ref{Example1} with $M=250$,  $N=1000$, $b=4$ shown on a wider domain $[-4,4]$ (two periods are shown).} 
\label{F3}
\end{figure}

\begin{figure}[h] \begin{center}
\includegraphics[width=2.4in,angle=0]{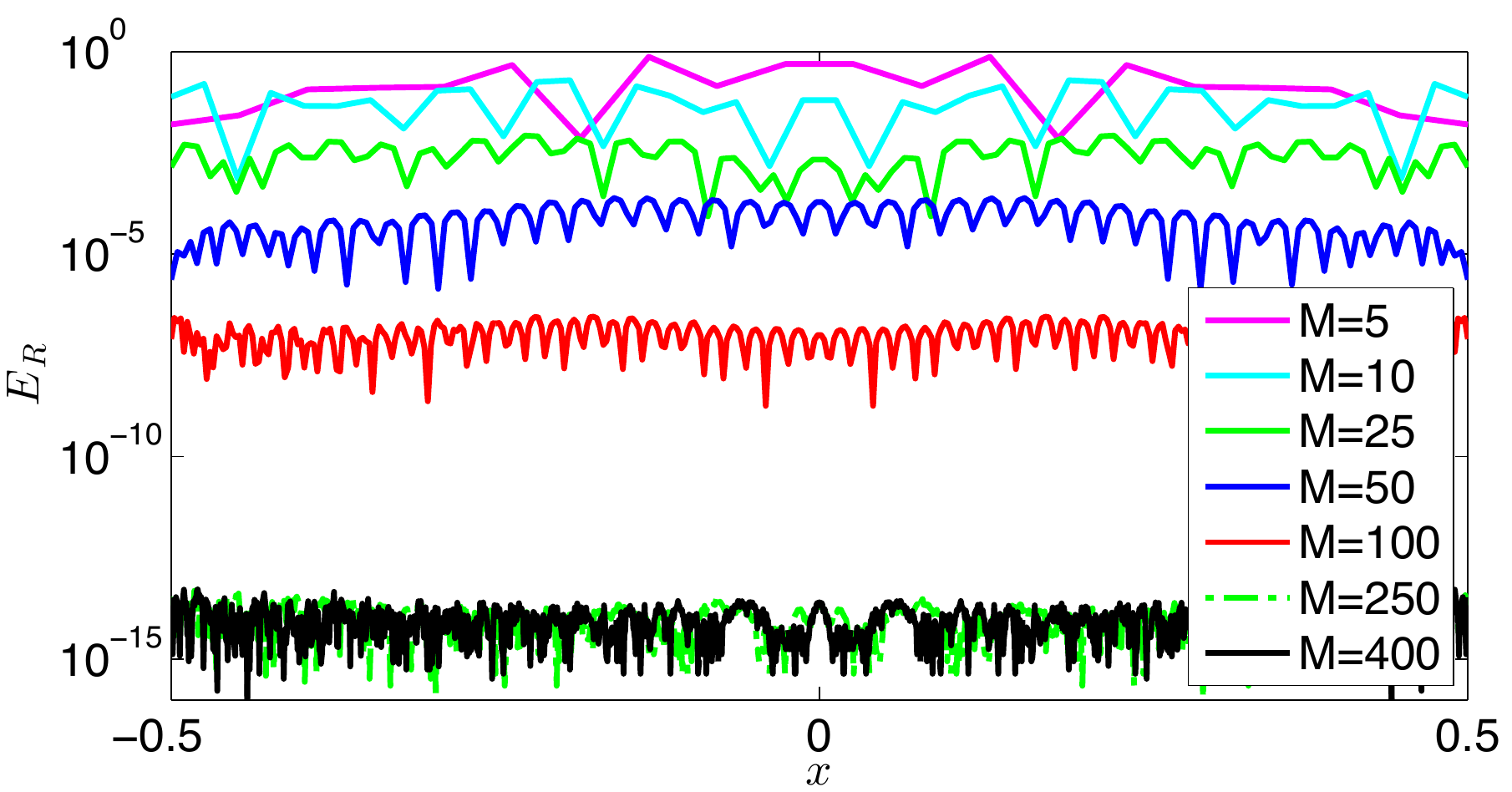}
\end{center}
\caption{Semilogy plot of the errors $E_R(x)$  in example \ref{Example1} on $[-0.5, 0.5]$ with $M=5$, $10$, $25$,  $50$, $100$ and $250$ and $400$, $N=4M$, $b=4$.} 
\label{F2}
\end{figure}

The two-pole transfer function with no delay \cite{Knockaert_Dhaene_2008}  is defined by
\begin{equation} \label{two_pole_eqn}
H(w)=\frac{r}{iw+s}+\frac{\overline{r}}{iw+\overline{s}}
\end{equation}
where $r=1+3i$, $s=1+2i$. Since both poles of this function located at $\pm 2+i$ are in the upper half plane,  the transfer function is causal as a linear combination of causal transforms. 
We sample data on the interval from $w=0$ to $w_{max}=6$, use the spectrum symmetry to obtain data on $[-w_{max},0)$ and scale the frequency interval from $[-w_{max}, w_{max}]$ to $[-0.5, 0.5]$. 
The real and imaginary parts of $H(x)$ 
are shown in Fig. \ref{F4}. Superimposed  are their causal Fourier continuations obtained using $M=250$, $N=1000$,  $b=4$ and solving the complex system (\ref{E3}) and its real counterpart (\ref{E4}). As can be seen, there is no essential difference in using complex or real formulation, though the real formulation (\ref{E4}) is slightly more ill-conditioned than the complex one. The data and the causal Fourier continuations are practically undistinguishable on $[-0.5, 0.5]$.  
%
To demonstrate the nature of continuations, we plot the same curves as in Fig. \ref{F4}  (only the real parts are presented, the imaginary parts have similar features) but on the extended domain $[-4,4]$ where we show two periods. These are plotted in Fig. \ref{F3}. It is obvious that the continuations oscillate in the extended region outside $[-0.5, 0.5]$. The frequency of these oscillations  increases with $M$. At the same time, the Fourier series become more and more accurate in approximation in the original interval $[-0.5, 0.5]$. To demonstrate this, we show the reconstruction errors $E_R(x)$,  defined in (\ref{error_Re}),  in Fig. \ref{F2} (semilogy plot is shown) in $[-0.5, 0.5]$ for various values of $M$ with $N=4M$ obtained using real formulation (\ref{E4}). The results for $E_I(x)$ are similar.
As $M$ increases from  $M=5$  to $250$, the order of the error decreases from  $10^{-1}$ to $10^{-14}$ for both real and imaginary parts. 
For example, with $M=250$, $N=1000$, $b=4$, both errors $E_R(x)$ and $E_I(x)$ are at the order of $4\times 10^{-14}$.
%
The results indicate that the error is uniform on the entire interval $[-0.5, 0.5]$ and does not exhibit boundary artifacts.
%
%
The above results demonstrate that the proposed technique is capable of verifying that the given data are causal with high accuracy. In this case, causality is satisfied with the error less than $10^{-13}$. These results are in agreement with the error estimates  (\ref{combined_bound_noise}),  (\ref{lambda12}) developed in Section \ref{error_analysis}.
 Since the data do not have any noise except of round-off errors and the transfer function is smooth on $[-0.5, 0,5]$, the reconstruction errors are dominated  by the fast decaying error in approximation of the smooth transfer function with its causal Fourier series for smaller $M$ and then by the error due to truncation of the singular values for high values of $M$, which it is close to the cut-off value $\xi=10^{-13}$. 
%
%
%
We also make another observation about the behavior of the errors $E_R(x)$ and $E_I(x)$ as $M$ increases for a causal smooth function. 
Even for small values of $M$, as $M$ doubles, the errors decrease by several orders of magnitude until the errors level off around $5\times 10^{-14}$ (see Table \ref{Tcausalerror}). This is a consequence of the fast convergence of a Fourier series for a smooth function.
\begin{table}[h]
\begin{center}
\begin{tabular}{|c|c|c||c|c|c|}
\hline
\rule{0cm}{10pt}
$N$ & $M$ & $||E_R||_\infty$, $||E_I||_\infty$ & $N$ & $M$ & $||E_R||_\infty$, $||E_I||_\infty$  \\[3pt]
\hline
\rule{0cm}{10pt}
$40$ & $10$ & $\sim 10^{-1}$ &  $400$ & $100$ & $\sim 10^{-7}$  \\[3pt]
\hline
\rule{0cm}{10pt}
$100$ & $25$ & $\sim 4\times 10^{-3}$ & $800$ & $200$  &$\sim 2\times 10^{-13}$ \\[3pt]
\hline
\rule{0cm}{10pt}
$200$ & $50$ & $\sim 10^{-4}$ & $1000$ & $250$ & $\sim  5\times10^{-14}$ \\[3pt]
\hline
\end{tabular}
\end{center}
\caption{The orders of errors  $E_R(x)$ and $E_I(x)$ in example \ref{Example1}  to demonstrate how fast reconstruction errors decay as resolution increases in case of a causal smooth function.}
\label{Tcausalerror}
\end{table} 

\begin{figure}[h] \begin{center}
\includegraphics[width=2.4in,angle=0]{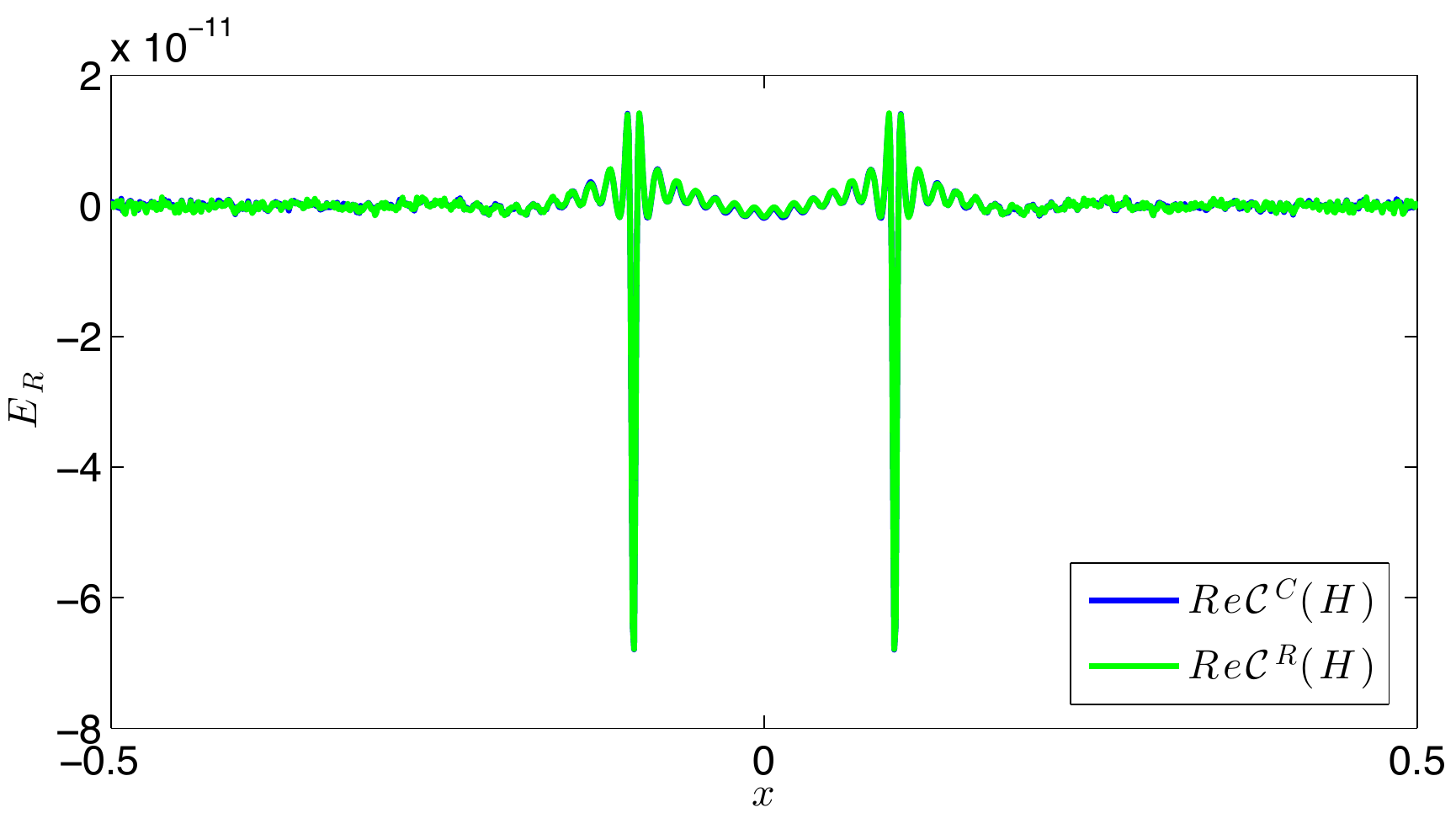}
\includegraphics[width=2.4in,angle=0]{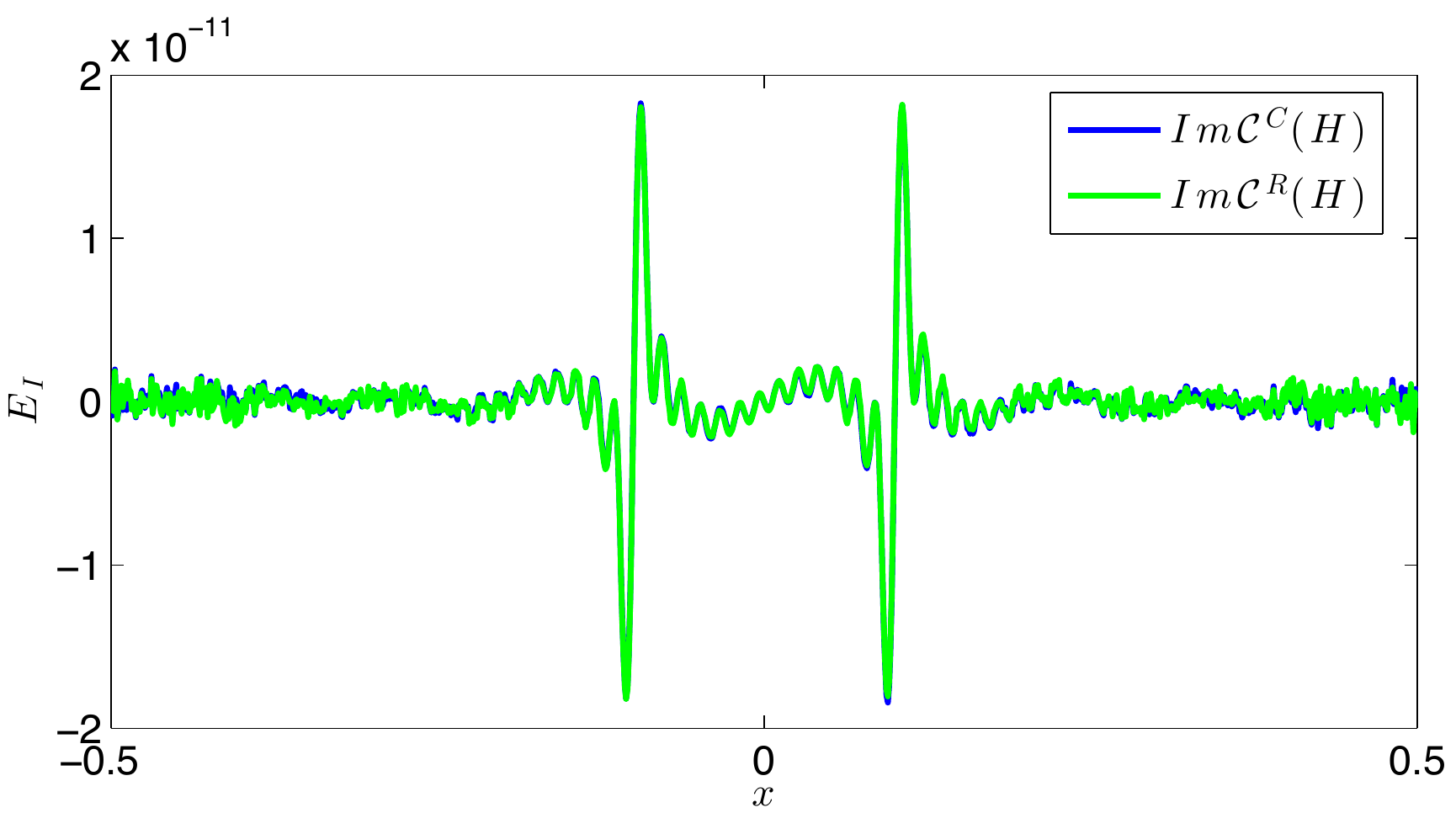}
\end{center}
\caption{Semilogy plot of $E_R(x)$ and $E_I(x)$  in example \ref{Example1} with Gaussian perturbation (\ref{Gauss_pert}) with $a=10^{-10}$ on  $[-0.5, 0.5]$ with $M=250$, $N=1000$, $b=4$.} 
\label{F5}
\end{figure}

Next we test how sensitive this method is to causality violations. We do this by imposing a localized non-causal perturbation on causal data. This artificial causality violation is modeled by a Gaussian function
\begin{equation} \label{Gauss_pert}
a\exp\left(-\frac{(x-x_0)^2}{2\sigma^2}\right) 
\end{equation}
of small amplitude $a$, centered at $x_0$ and added to $\Re H$, while keeping $\Im H$ unchanged. This type of non-causal perturbation was used in \cite{Triverio_Grivet_Talocia_2006} to test causality verification technique based on the generalized dispersion relations.  We use the Gaussian centered at $x_0=0.1$ with standard deviation $\sigma=10^{-2}/6$, so its ``support" is concentrated on a very narrow interval of length $10^{-2}$ and outside this interval the values of the perturbation are very close to $0$. The advantage of using a Gaussian perturbation is that it can be localized and allows one to verify if the proposed method is capable of detecting the location of causality violation. By varying the amplitude $a$, we can impose larger or smaller causality violations. The smaller $a$ can be used, the more sensitive method is for detecting causality violations.
With $a=10^{-10}$, which gives a very small causality violation, the error between the data and its causal Fourier continuation is shown in Fig. \ref{F5}. It is clear that the error has very pronounced spikes 
at $x=\pm 0.1$ due to symmetry that correspond to the exact locations of the Gaussian perturbation. These spikes are of the order $10^{-11}$   whereas on the rest of the interval  the error is about $10$ times smaller. For larger perturbations, the results are similar. For example, with $a=10^{-8}$, the error at $\pm 0.1$ is of the order of $10^{-9}$ and the rest of the interval has the error  $10$ times smaller, etc. We can see that the reconstruction error in this case is strongly dominated by the error (perturbation) in the data at the location of causality violation and the magnitude of the error is of the same order as the order of the perturbation. At the same time, the transfer function itself is very smooth, which ensures fast convergence of the Fourier series. The results are in a perfect agreement with the error estimate (\ref{combined_bound_noise}),  (\ref{lambda12}).
%
%

Another perturbation that we consider is a cosine function
%
$a\cos(20\pi x)$
%
that we also add to $\Re H(x)$ but keep $\Im H(x)$ unaltered. Adding a non-causal cosine perturbation makes the transfer function non-causal on the entire  interval and higher reconstruction errors are expected everywhere.
We find that both errors $E_R(x)$ and $E_I(x)$
oscillate with the frequency and amplitude similar to those of the perturbation. For example, with $a=10^{-10}$, these errors are of the order  $6\times 10^{-11}$.
Increasing/decreasing $a$ increases/decreases the magnitude of the error. 
%
%
To see how the level of a noise can be detected, set, for example, $a=10^{-5}$ and compute the reconstruction errors $E_R(x)$ and $E_I(x)$ as the resolution increases. We vary $M$ from $10$ to $250$ and analyze the order of the errors. The results are shown in Table \ref{two_pole_noncausal}. As we can see, the reconstruction errors first decrease as the resolution (the number $N$ of points) and the number $M$ of Fourier coefficients  increase until the error reaches the size of the perturbation, which happens at $M=100$, and after that the order of the error  does not decrease further and levels off instead.   
\begin{table}[h]
\begin{center}
\begin{tabular}{|c|c|c||c|c|c|}
\hline
\rule{0cm}{10pt}
$N$ & $M$ & $||E_R||_\infty$, $||E_I||_\infty$ & $N$ & $M$ & $||E_R||_\infty$, $||E_I||_\infty$  \\[3pt]
\hline
\rule{0cm}{10pt}
$40$ & $10$ & $\sim 10^{-1}$ & $400$ & $100$ & $\sim 7\times 10^{-6}$  \\[3pt]
\hline
\rule{0cm}{10pt}
$100$ & $25$ & $\sim 4\times 10^{-3}$ & $800$ & $200$  &$\sim 7\times 10^{-6}$ \\[3pt]
\hline
\rule{0cm}{10pt}
$200$ & $50$ & $\sim 10^{-4}$ & $1000$ & $250$ & $\sim  7\times10^{-6}$ \\[3pt]
\hline
\end{tabular}
\end{center}
\vskip5pt
\caption{The orders of errors  $E_R(x)$ and $E_I(x)$ in example \ref{Example1} with non-causal cosine perturbation 
$a\cos(20\pi x)$
with $a=10^{-5}$ as $M$ and $N$ increase in proportion $N=4M$.
%
}
\label{two_pole_noncausal}
\end{table} 
With a smaller perturbation amplitude $a$, the error levels off at a larger value of $M$. For example, with $a=10^{-10}$, the reconstruction errors stop decreasing at $M=200$.

\subsection{Finite Element Model of a DRAM package} \label{example_Micron}

\begin{figure}[h] \begin{center}
\includegraphics[width=2.4in,angle=0]{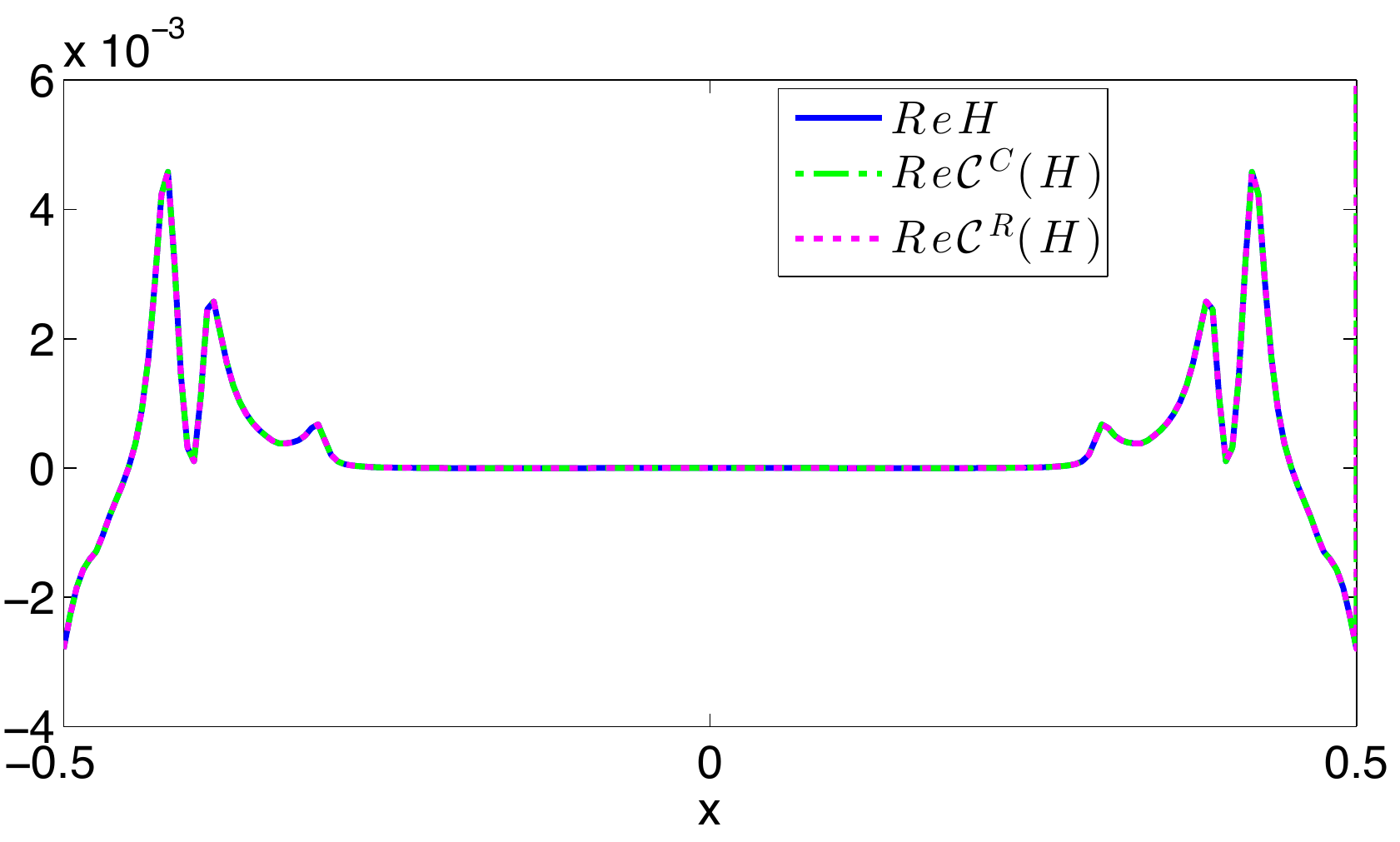}
\includegraphics[width=2.4in,angle=0]{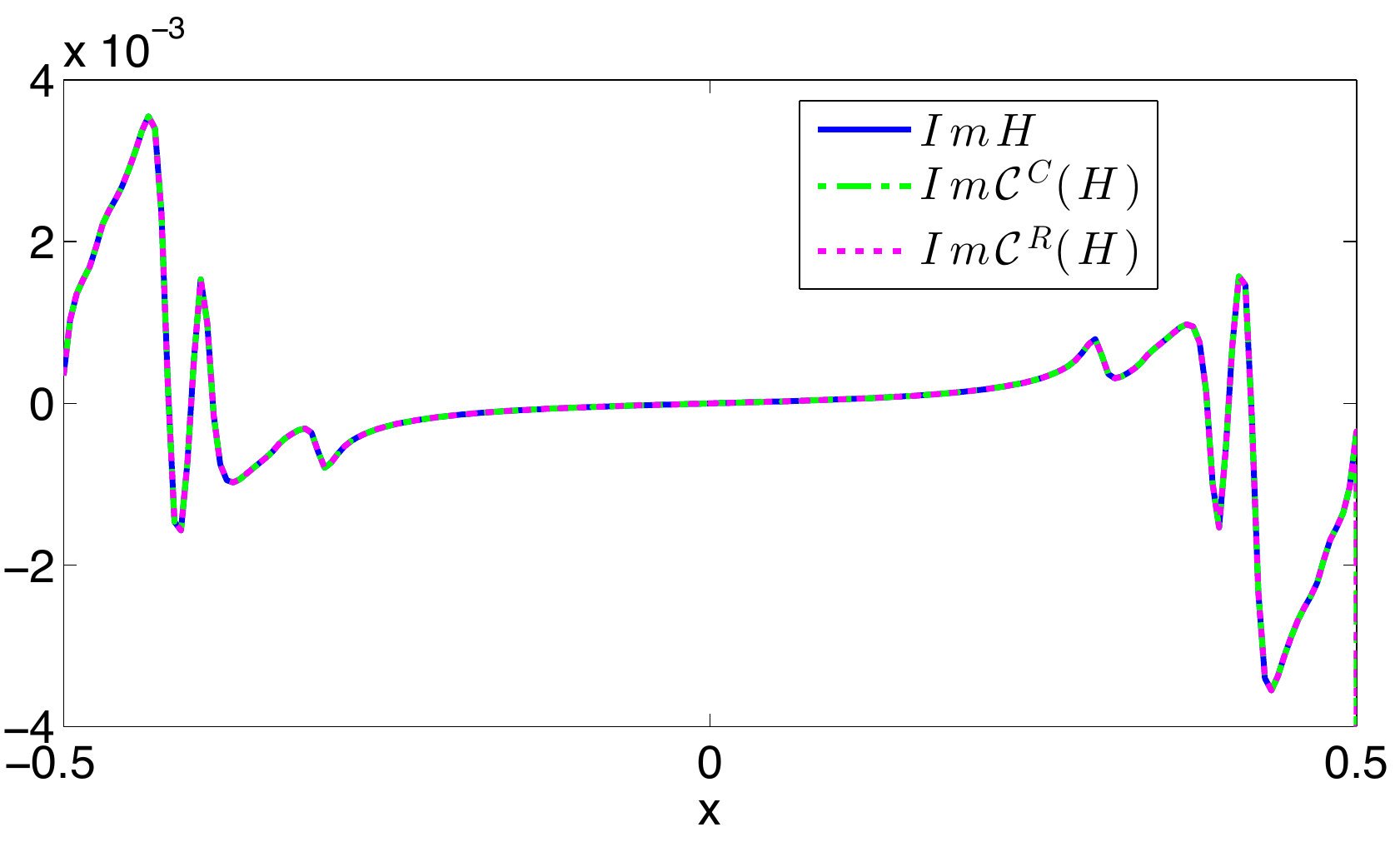}
\end{center}
\caption{$H(x)$ with its causal Fourier continuations ${\mathcal C}^C(H)$ and  ${\mathcal C}^R(H)$ in example \ref{example_Micron} with $M=100$,  $N=199$,  $b=1.1$.} 
\label{FMicron}
\end{figure}

In this example we use a scattering matrix $S$ generated by a Finite Element Modeling (FEM) of a DRAM package (courtesy of Micron). The package contains $110$ input and output ports. The simulation process was performed for $100$ equally spaced frequency points ranging from $w_{min}$=0 to $w_{max}=5$ GHz. We expect data to be causal but, perhaps, with some error due to limited accuracy of  numerical simulations.
For simplicity, we apply our method to the $S$-parameter $S(100,1)$ rather than the entire $110 \times 110$ $S$-matrix. The selected $S$-parameter $H(w)=S(100,1)$  relates the output signal from port $100$ to the input signal at port $1$ as a function of frequency $w$. Our approach can be extended to the entire $S$-matrix by applying the method to every element of the scattering matrix $S$. The graph of $H(x)$ is shown in Fig. \ref{FMicron}. 
In this test example, the number of samples in $[-0.5, 0.5]$ is fixed at $N=2\cdot 100-1=199$ (by symmetry), so to construct a causal Fourier continuation we use $M=100$ Fourier coefficients  and  can only vary the length $b$ of the extended region. 
Because the transfer function has oscillations and  high slopes in the boundary region, we have to use a smaller value of $b$. We find  $b=1.1$ to be  optimal  by trying a few different values of $b\in(1,2)$.
Slightly higher values of $b$ give similar results while using too large $b$ does not produce small enough error, most likely because we have fixed resolution and more data points would be needed to construct a causal continuation on a larger domain with good enough resolution. 

\begin{figure}[h] \begin{center}
\includegraphics[width=2.4in,angle=0]{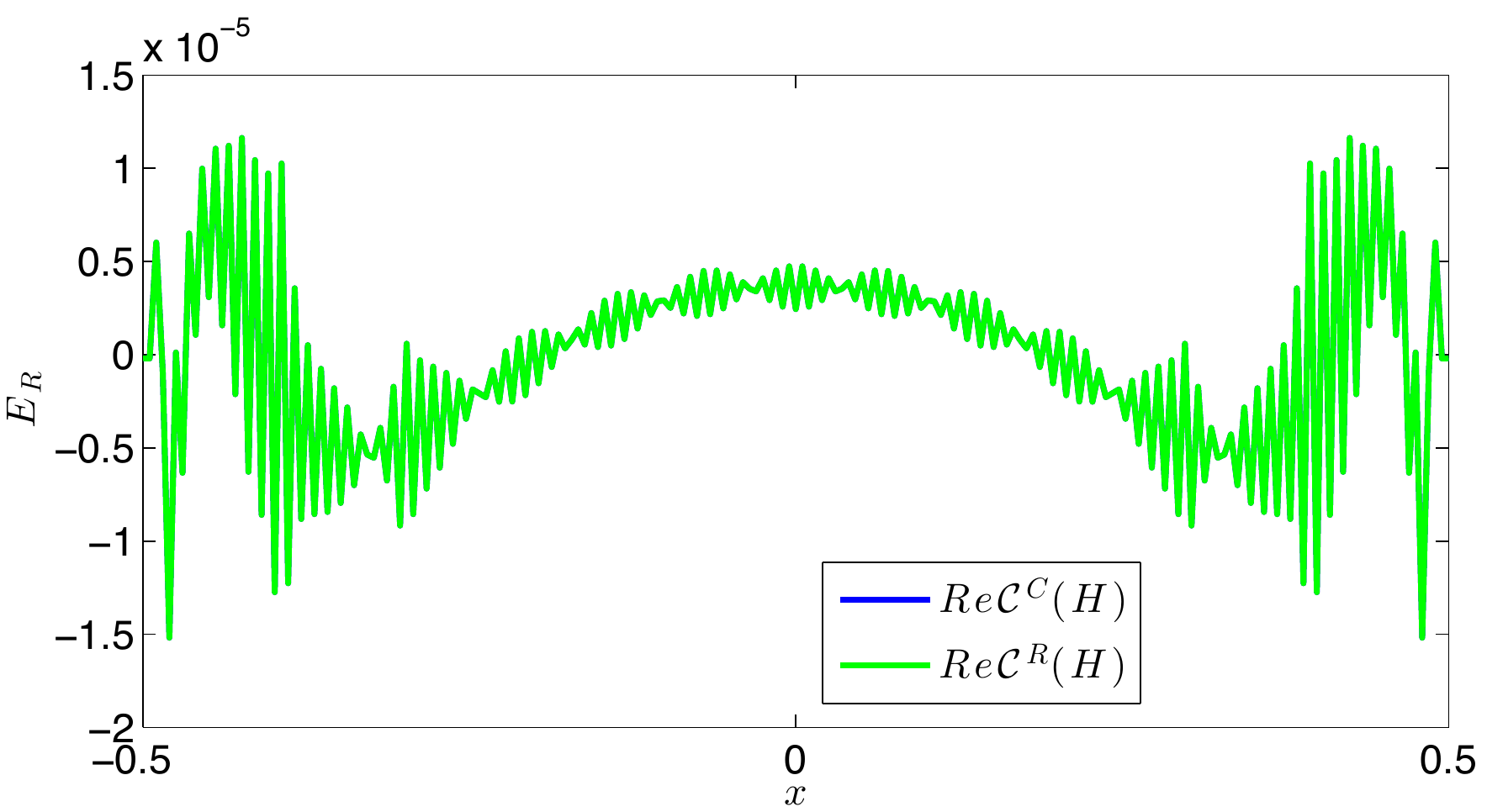}
\end{center}
\caption{Errors $E_R(x)$  in approximation of $S(100,1)$ in example \ref{example_Micron} with $N=199$, $M=100$  and $b=1.1$.  Errors $E_I(x)$ have the same order.
} 
\label{FMicron_errors}
\end{figure}

The  errors $E_R(x)$  with the above parameters are shown in Fig. \ref{FMicron_errors}. The errors $E_I(x)$ have the same order of $10^{-5}$. This  
implies that the dispersion relations are satisfied within error  on the order of $10^{-5}$.  Since the data came from finite element simulations, we expect their accuracy to be on the order of $10^{-6}$ or $10^{-7}$ at least. To verify if the relatively large reconstruction errors comes from a noise in the data or a causal Fourier series approximation error, we use data with $N$ samples, then every other and every forth samples, i.e. with $N$, $\frac{N-1}{2}+1$ and $\frac{N-1}{4}+1$ samples to find a model (\ref{error_epsilon1}) of the form $CM^\alpha$ using the least squares method. We find that both $E_R$ and $E_I$ decay approximately as ${\mathcal O}(M^{-3})$.
These models were extrapolated to higher values of $M$ as shown in Fig. \ref{FMicron_errors}, where we plot the $l_2$ norms of actual reconstruction errors and fitted model curves. The extrapolated error curves indicate that the error may be decreased further if more frequency samples are available. In this case, we say that the transfer function $H(x)$ satisfies dispersion relations within $10^{-5}$, i.e. the transfer function $H(x)$ is causal within the error at most $10^{-5}$, and the causality violations, if present, are smaller than or at most on the order of $10^{-5}$. To determine the actual level of a noise in the data, higher resolution of frequency responses would be needed.

For comparison, using periodic polynomial continuation method \cite{Aboutaleb_Barannyk_Elshabini_Barlow_WMED13, Barannyk_Aboutaleb_Elshabini_Barlow_IMAPS} with $8$th degree polynomial applied to the transfer function in this example, the error in approximation of $H(x)$ is about $2\times 10^{-3}$ that is by two orders of magnitude larger than with the spectral continuation, it is not uniform and the largest in the boundary domain. 
\begin{figure}[h] \begin{center}
\includegraphics[width=2.4in,angle=0]{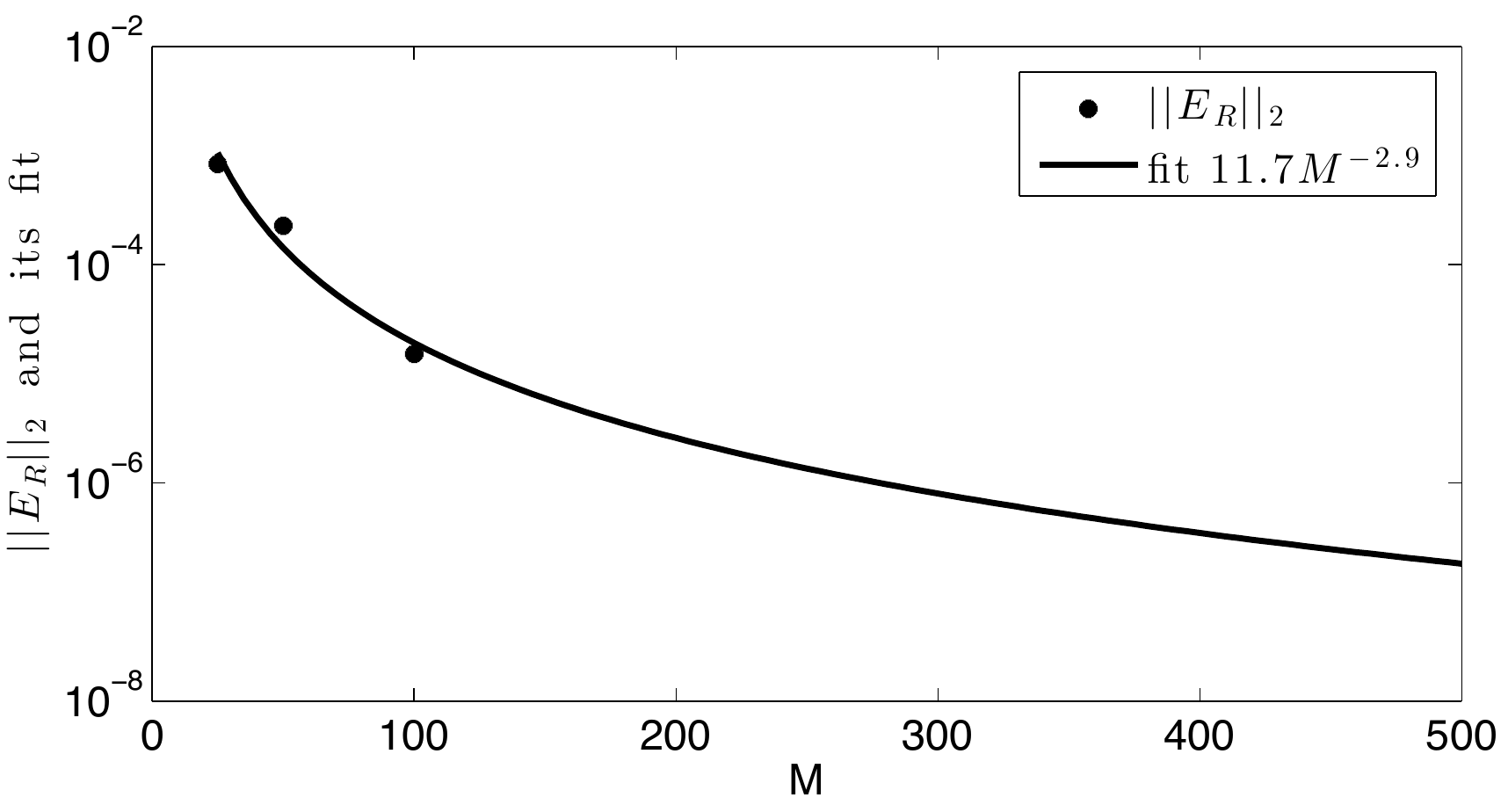}
\end{center}
\caption{Errors $||E_R(x)||_2$  in approximation of $S(100,1)$ in example \ref{example_Micron} with $N=199$, $99$, $49$ and $M=100$, $50$, $25$, respectively,  and $b=1.1$, together with their least squares fit: $||E_R||_2\sim 11.7M^{-2.9}$ and extrapolation for higher values of $M$. For $||E_I(x)||_2$ we find $||E_I||_2\sim 23.5M^{-3.1}$.
} 
\label{FMicron_errors_fit}
\end{figure}

\subsection{Transmission line example} \label{transmission_line_example}

\begin{figure}[h] \begin{center}
\includegraphics[width=2.4in,angle=0]{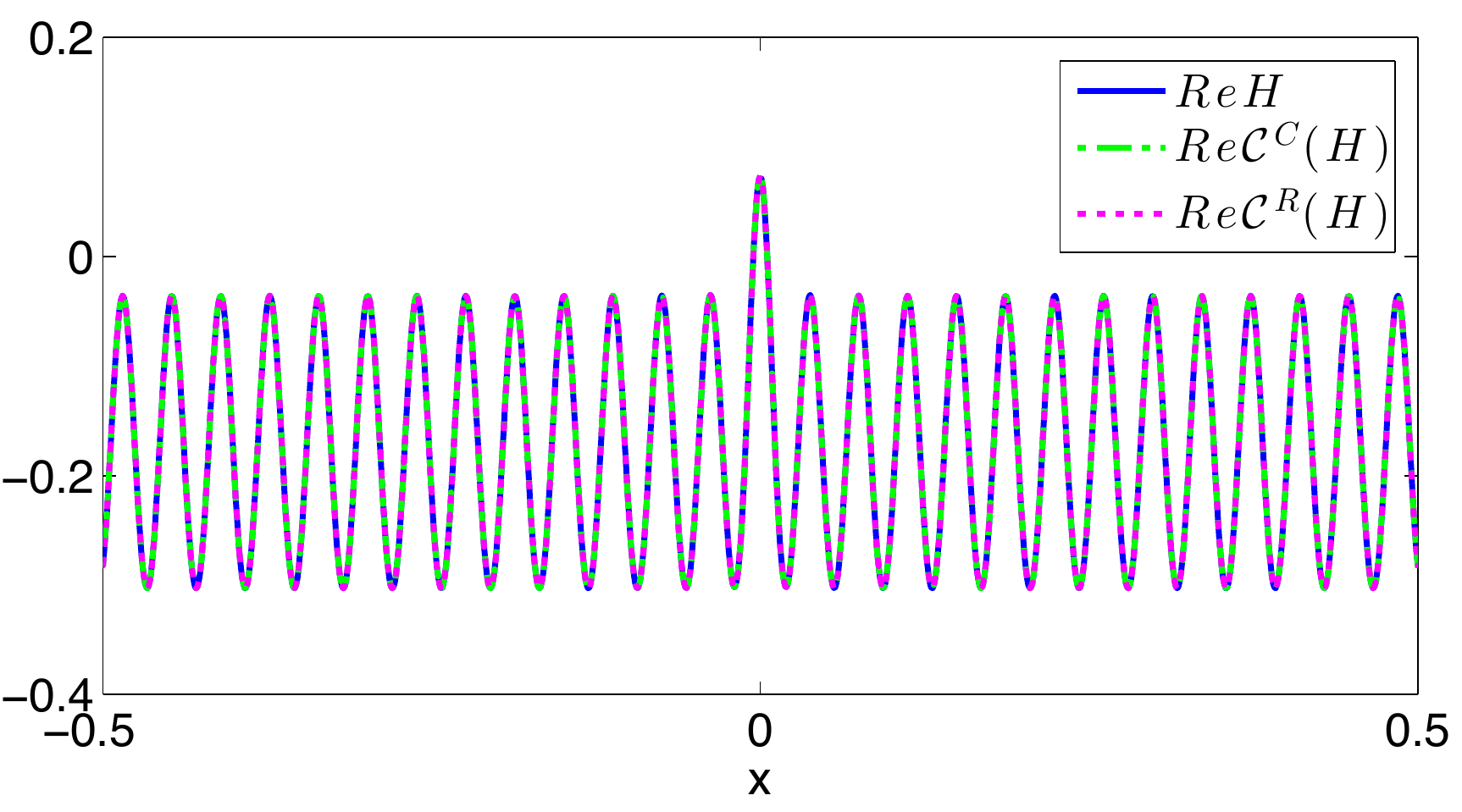}
\end{center}
\caption{Transfer function $H(x)$ and its causal Fourier continuations ${\mathcal C}^C(H)$ and  ${\mathcal C}^R(H)$  in example \ref{transmission_line_example} with $M=1500$, $N=3000$,  $b=4$ (real parts are shown).} 
\label{F9}
\end{figure}

We consider a uniform transmission line segment that has the following per-unit-length parameters: $L=4.73$ nH/cm, $C=3.8$ pF/cm, $R=0.8$ $\Omega$/cm, $G=0$ and length ${\mathcal L}=10$ cm. The frequency is sampled on the interval $(0, 5.0]$ GHz. This example was used in \cite{Triverio_Grivet_Talocia_2006} to analyze causality using generalized dispersion relations with subtractions. The scattering matrix of the structure was computed using Matlab function {\tt rlgc2s}. Then we consider the element $H(w)=S_{11}(w)$. Due to limitation of the model used in the  function 
{\tt rlgc2s}, we cannot obtain the value of the transfer function at $w=0$  (DC) but we can sample it from any small nonzero frequency. It is typical for systems to have frequency response missing at low frequencies and can occur either during measurements or simulations. However, the value of $H(w)$ at $w=0$ is finite, because the magnitude of $S_{11}$ must be bounded by $1$. Hence, we  have a bandpass case. Once we choose the number of points and the corresponding $w_{min}>0$, we can sample frequencies from $[w_{min}, w_{max}]$. We use $w_{max}=5.0$ GHz. Using symmetry conditions we reflect the values of the transfer function for negative frequencies as for the baseband case considered above. We know that $\Im H(w)$ equals $0$ at $w=0$ but $\Re H(w)$ is to be computed.  Since we do not have a value of $\Re H(w)$ at $w=0$, our frequencies at which the values of the transfer function are available will have a gap at $w=0$. Nevertheless, our approach is still applicable since it does not require the data points to be equally spaced. In this example,  we get better results, however, with smaller $w_{min}$. Alternatively, we can use a polynomial interpolation to find a value of the real part of $H(w)$ at $w=0$. The value of the imaginary part $\Im H(0)=0$ by symmetry. This approach is not very accurate since it does not take into account causality when the polynomial interpolation of $\Re H(w)$ is constructed, and produces a larger error compared with just skipping the value at $w=0$ and using spectral continuation approach directly. With our technique and $M=1500$, $N=3000$,  $b=4$ we are able to construct a causal Fourier continuation accurate within $3\times 10^{-15}$. The graphs of $\Re H(w)$ 
together with their causal Fourier continuations  are presented in Fig. \ref{F9}. Agreement for $\Im H(w)$  is similar.
With smaller values of $w_{max}$, it is enough to use smaller values of $M$ and $N$ in the same proportion $2M=N$ to get the same order of accuracy. For example, to get an error in approximation of the transfer function on the original interval at the order of $10^{-14}$ and $w_{max}=3.0$ GHz, it is enough to use $M=750$, $N=1500$, while for $w_{max}=1$ GHz one could use $M=250$, $N=500$. In this case, just more data points would be needed to have the same resolution on the longer domain.
%
%
Computations presented in this paper were done using Matlab on a MacBook Pro with 2.93\,GHz Intel Core 2 Duo Processor and 4 GB RAM. The CPU times (in seconds) for computing SVD, minimum norm solution of a linear system using the least squares method as well as overall CPU time for constructing a causal Fourier continuation with $M$ varying from $50$ to $1500$ with $N=2M$ are shown in Table  \ref{T_CPU}. 
%
\begin{table}[h]
\begin{center}
\begin{tabular}{|c|c|c|c|c|}
\hline
\rule{0cm}{10pt}
$N$ & $M$ &      CPU time   &  CPU time   & CPU time   \\[3pt]
& & SVD & minnmsvd & continuation \\[3pt]
\hline
\rule{0cm}{10pt}
100 & 50   &  0.0091     &     0.0141        &      0.0299 \\[3pt]
\hline
\rule{0cm}{10pt}
200  & 100 &  0.1910     &     0.0340        &      0.0781 \\[3pt]
\hline
\rule{0cm}{10pt}
500   & 250 & 0.9626     &     0.1646       &       0.4040 \\[3pt]
\hline
\rule{0cm}{10pt}
1000  & 500 &  3.9779     &     1.5061       &       3.9063 \\[3pt]
\hline
\rule{0cm}{10pt}
2000  & 1000&  25.5921   &      9.8634      &        35.9250 \\[3pt]
\hline
\rule{0cm}{10pt}
3000  & 1500 &  68.6258    &     48.4047   &  117.6096  \\[3pt]
\hline
\end{tabular}
\end{center}
\vskip5pt
\caption{CPU times (in seconds) for computing SVD, finding minimum norm solution using the least squares method and overall CPU time for constructing causal Fourier continuation for $M$ varying from $50$ to $1500$ with $N=2M$.}
\label{T_CPU}
\end{table} 


%
\begin{figure}[h] \begin{center}
\includegraphics[width=2.4in,angle=0]{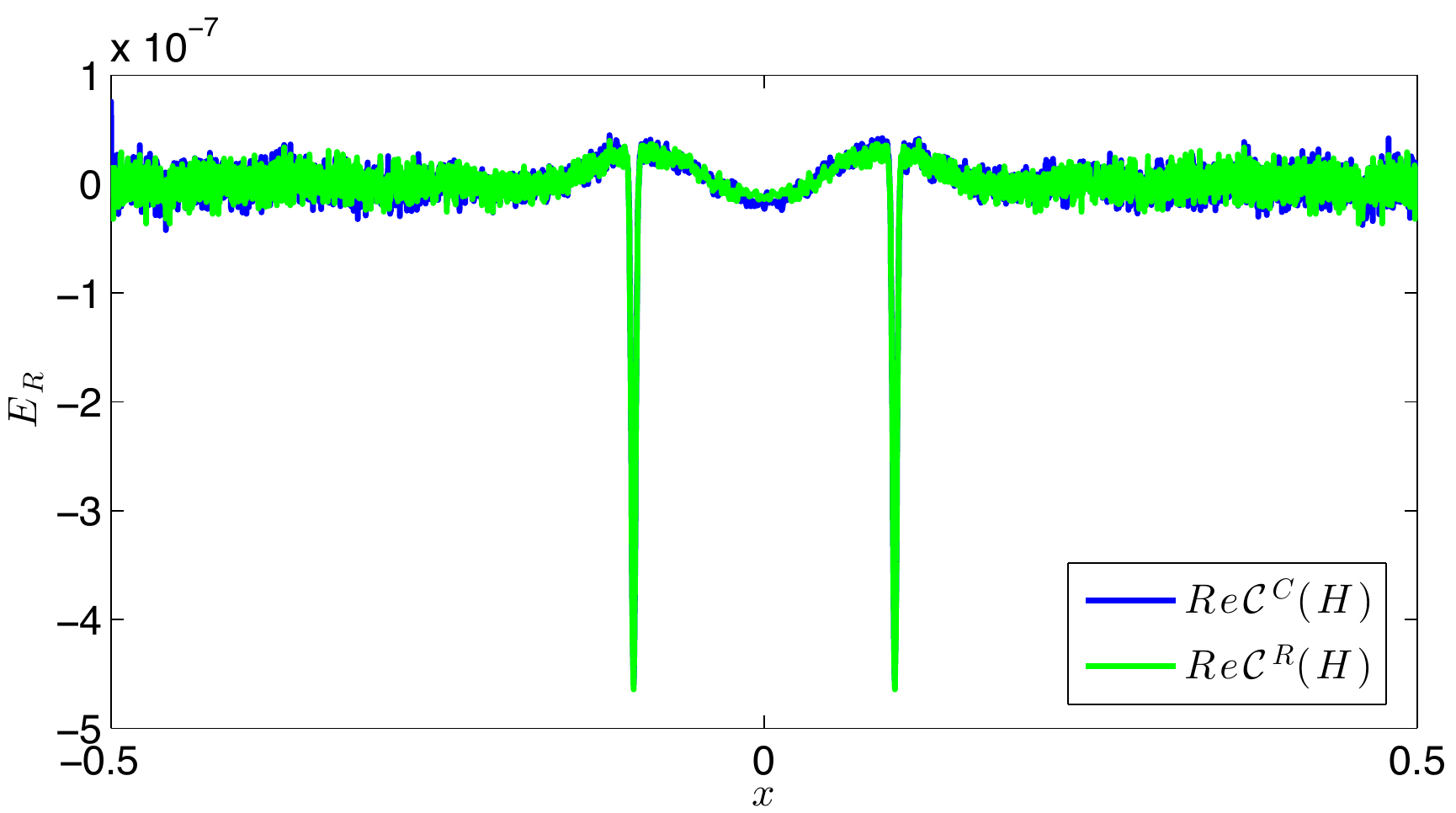}
\end{center}
\caption{Errors $E_R(x)$ 
in example \ref{transmission_line_example} with $M=1500$, $N=3000$, $b=4$ and non-causal Gaussian perturbation (\ref{Gauss_pert}) with  $a=10^{-6}$, $\sigma=10^{-2}/6$, centered at $x_0=0.1$. The errors $E_I(x)$   have similar spikes at causality violation locations.} 
\label{F13}
\end{figure}

Now we impose a Gaussian perturbation (\ref{Gauss_pert}) to the transfer function  as in example \ref{Example1} with amplitude $a=10^{-6}$ to model a causality violation  located at $x_0=0.1$ 
with the standard deviation $\sigma=10^{-2}/6$, so that the ``support" of the Gaussian is approximately $10^{-2}$ or $1/10$ of the entire bandwidth. This results in the error having very pronounced spikes, shown in Fig. \ref{F13}, of magnitude $4.5\times 10^{-7}$  in both $\Re H$ and $\Im H $ at the locations of the perturbation, while as before the error on the rest of the interval is approximately $10$ times smaller.
%
%
We find that we are able to detect reliably a causality violation of amplitude  up to $a=5\times 10^{-14}$. 

 
In papers \cite{Triverio_Grivet_Talocia_2006, Triverio_Grivet_Talocia_2006_2, Triverio_Grivet_Talocia_2008} an excellent work is done by utilizing the generalized dispersion relations with subtractions to check causality of raw frequency responses. The performed error analysis provides explicit frequency dependent error estimates to account for finite frequency resolution (discretization error in computing Hilbert transform integral) and finite bandwidth (truncation error due to using only a finite frequency interval instead of the entire frequency axis) and unbiases the causality violations from numerical discretization and domain truncation errors. It is shown that by using more subtractions, one can make the truncation error arbitrary small but the discretization error does not go away since it is fixed by the resolution of given frequency responses. The authors report that if a causality violations are too small (with amplitude smaller than $10^{-5}$) and smooth, using more subtractions may not affect the overall error since it is then dominated  by the discretization error. In addition, even with placing more subtraction points in the boundary regions close to $\pm w_{max}$, the truncation error 
always diverges at the bandwidth edges due to the missing out-of-band samples.
%
In the current work, we are able to remove boundary artifacts and  detect  really small localized infinitely smooth (Gaussian) causality violations. The resolution of data is also important since the number $N$ of collocation points at which frequency responses are available, dictates the number $M$ of Fourier coefficients with $N=2M$, and, hence, the sensitivity of the method to causality violations. 
For instance, in the above example that was also used in \cite{Triverio_Grivet_Talocia_2006} (but we use a smaller amplitude $a=10^{-6}$ and the Gaussian with more narrow $6\sigma=10^{-2}$ ``support"), with $N=250$ ($125$ points in $[0,0.5]$) we can detect causality violation with $a=5\times 10^{-7}$, whereas  $N=500$ and $N=1000$ are capable of detecting violations with $a=10^{-12}$ and $a=5\times 10^{-13}$, respectively. Similarly to \cite{Triverio_Grivet_Talocia_2008}, we also find that it is more difficult to detect wide causality violations. As the ``support" $6\sigma$ of the Gaussian  increases, the spikes in the reconstruction errors become wider and shorter, and eventually it is not possible to determine the location of the causality violation since the error at causality violation locations has the same order as the error on the rest of the interval. For $6\sigma\leq 0.1$,  the spikes in the errors are still observable, whereas for a bigger value $6\sigma \geq 0.2$, the reconstruction errors are uniform. Increasing the resolution of the data does not decrease the reconstruction errors, that indicates the presence of causality violations of the order of reconstruction errors.

\subsection{Delayed Gaussian example} \label{delayed_Gaussian}

\begin{figure}[h] \begin{center}
\includegraphics[width=2.4in,angle=0]{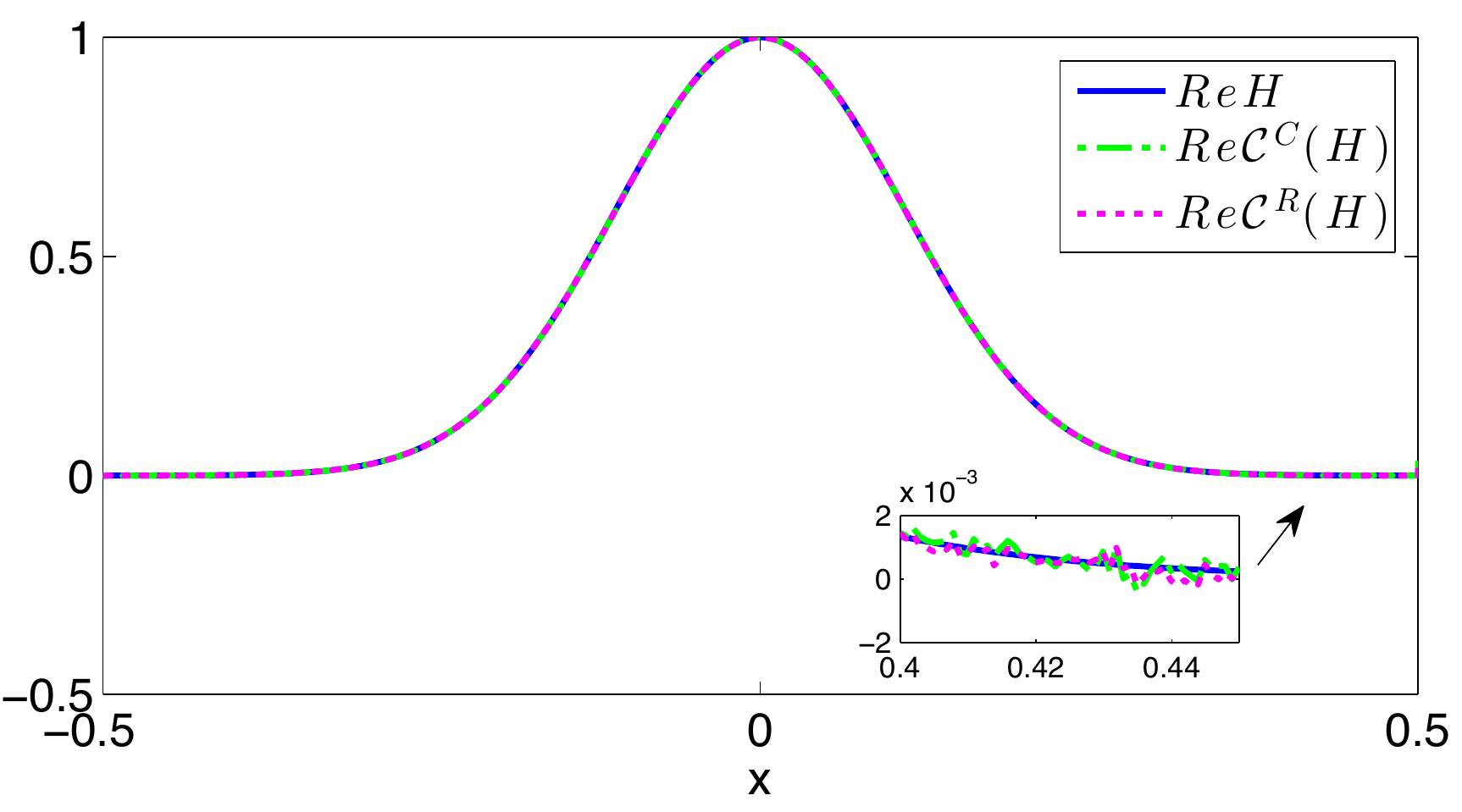}
\end{center}
\caption{Noncausal transfer function $H(w,t_d)$ in example \ref{delayed_Gaussian} and its Fourier continuations ${\mathcal C}^C(H)$ and  ${\mathcal C}^R(H)$  with $M=250$, $N=500$,  $b=2$, $t_d=0.1\sigma$ (only real parts are shown).} 
\label{F15}
\end{figure}

\begin{figure}[h] \begin{center}
\includegraphics[width=2.4in,angle=0]{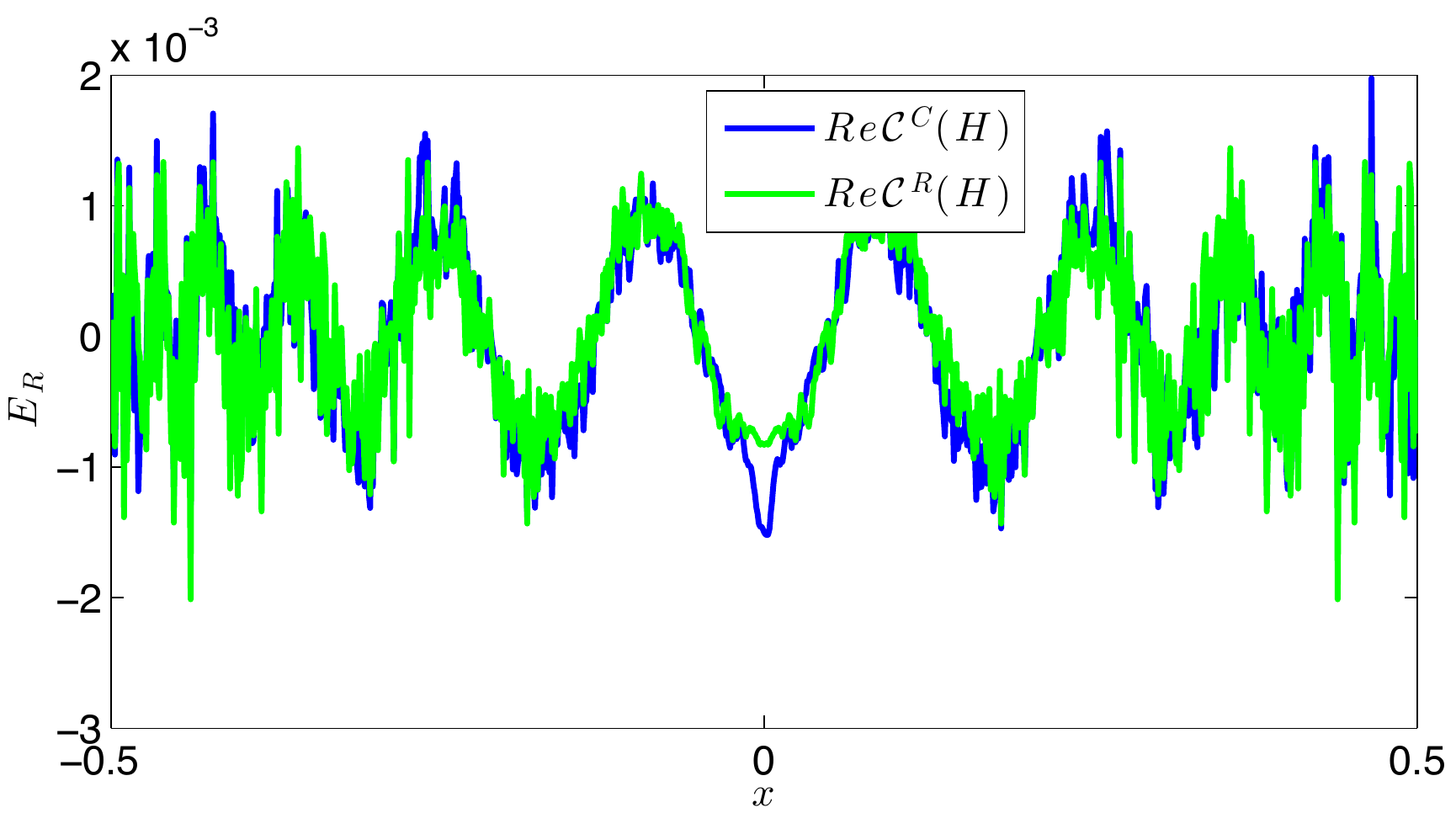}
\end{center}
\caption{Errors $E_R(x)$  between the noncausal transfer function $H(w,t_d)$  and its causal Fourier continuations ${\mathcal C}^C(H)$ and  ${\mathcal C}^R(H)$ in example \ref{delayed_Gaussian} with $M=250$, $N=500$,  $b=2$, $t_d=0.1\sigma$. Errors  $E_I(x)$ have the same order.} 
\label{F16}
\end{figure}

Here we test the performance of the method with an example of a delayed Gaussian function that was used in \cite{Xu_Zeng_He_Han_2006} to check causality of interconnects through the minimum-phase and all-pass decomposition. We consider the impulse response function modeled by a Gaussian with the center of the peak $t_d$ and standard deviation $\sigma$:
\[
h(t,t_d)=\exp\left[-\frac{(t-t_d)^2}{2\sigma^2}\right].
\]
If $t_d=0$, the Gaussian function $h(t,0)$ is even, so it cannot be causal. As $t_d$ increases, the center of the peak moves to the right 
and for $t_d>3\sigma$ the impulse response function $h(t,t_d)$ can be gradually made causal.  
%
The corresponding transfer function is
\[
H(w,t_d)=\exp\left[-2(\pi w \sigma)^2-2i\pi w\, t_d\right]
\]
which is a periodic function damped by an exponentially decaying function. We consider two regimes. One has value of $t_d<3\sigma$, so that the transfer function $H(w,t_d)$ is non-causal. In the second regime, the delay $t_d> 3\sigma$ is big enough   to make  the transfer function $H(w,t_d)$ causal. We fix $b=2$, $\sigma=2$ and sample $w$ from the interval $[0,3.6\times 10^{8}]$  Hz  and consider first the case with $t_d=0.1\sigma$. The real  part of $H(w,t_d)$ is shown in Fig. \ref{F15} together with its Fourier continuations. One can clearly see that Fourier continuations do not match well $\Re H$. Instead, there are visible high frequency oscillations throughout the domain
as confirmed by analyzing  the reconstruction errors $E_R(x)$. With  $M=250$, $N=500$, the magnitude of the errors is about $2\times 10^{-3}$  (see Fig. \ref{F16}). When $N$ is increased in proportion $N=2M$, the error slightly increases. For example, with $M=1000$, $N=2000$, the errors are about $4\times 10^{-3}$.
%
%
Varying the length $b$ of the extended domain does not decrease the magnitude of the error. This is in agreement with the error estimate (\ref{combined_bound_noise}), (\ref{lambda12})  since in this case the reconstruction error  is dominated by causality violations. Results for $E_I(x)$ are similar. 

\begin{figure}[h] \begin{center}
\includegraphics[width=2.4in,angle=0]{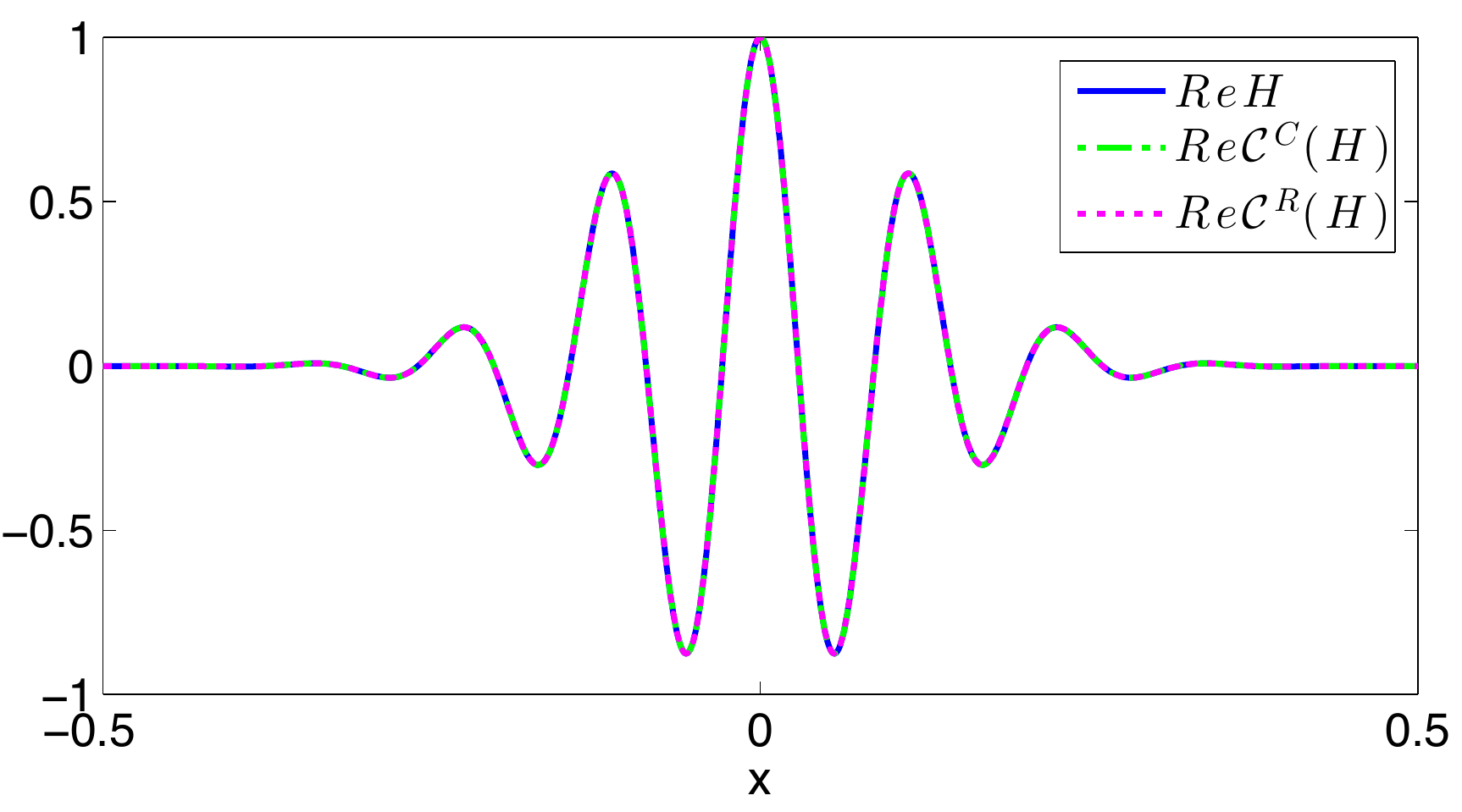}
\end{center}
\caption{Causal transfer function $H(w,t_d)$ in example \ref{delayed_Gaussian} and its Fourier continuation with $M=250$, $N=500$,  $b=2$, $t_d=6\sigma$ (real parts are shown).} 
\label{F17}
\end{figure}

\begin{figure}[h] \begin{center}
\includegraphics[width=2.4in,angle=0]{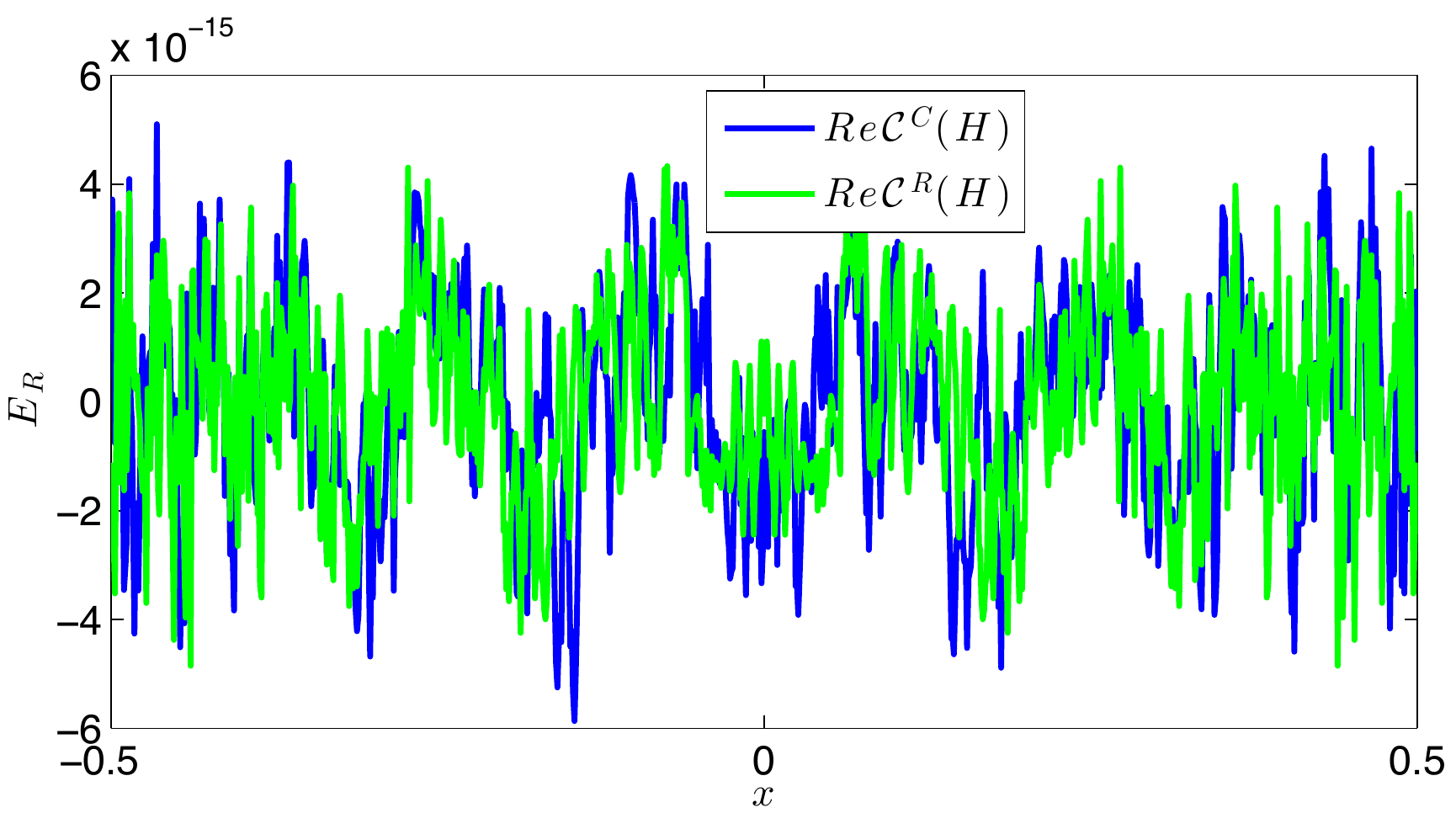}
\end{center}
\caption{Errors $E_R(x)$ between the causal transfer function $H(w,t_d)$  and its causal Fourier continuation in example \ref{delayed_Gaussian} with $M=250$, $N=500$,  $b=2$, $t_d=6\sigma$. Errors $E_I(x)$ have the same order.} 
\label{F18}
\end{figure}

In the second case, we set $t_d=6\sigma$, which should give a causal transfer function. $\Re H(w,t_d)$  together with its Fourier continuations are shown in Fig.  \ref{F17}.  Both reconstruction errors $E_R(x)$ and $E_I(x)$ drop to the order of $2\times 10^{-15}$ (see Fig. \ref{F18}).  
In this case, the transfer function is causal and infinitely smooth, so the error in approximation of such function with a causal Fourier series decays quickly even for moderate values of $M$. 

We do observe a gradual change of the  non-causal Gaussian function into a causal function as $t_d$ increases. Writing $t_d=\gamma\sigma$ and vary  values $\gamma=1$, $2$, $4$, and $5$ and find that the $l_\infty$ norms of both errors $E_R$ and $E_I$ are $5\times 10^{-6}$, $10^{-7}$, $3\times 10^{-12}$ and $10^{-14}$, respectively, i.e. the error decays as $t_d$ increases as  expected. 
%
%

\section{Conclusions} \label{conclusions}

We present a numerical method for verification and enforcement, if necessary, of causality of bandlimited tabulated frequency responses, that can be employed before the data are used for macromodeling. The approach is based on  the Kramers-Kr\"onig  dispersion relations and a construction of SVD-based causal Fourier continuations. This is done by calculating accurate causal Fourier series approximations of transfer functions, not periodic in general, and allowing the causal Fourier series to be periodic in an extended domain.
%
%
The causality is imposed directly on Fourier coefficients using dispersion relations that require real and imaginary parts of a causal function to be a Hilbert transform pair. The approach eliminates the necessity of approximating the behavior of the transfer function at infinity, which is known to be a source of significant errors in computation of the Hilbert transform defined on an infinite domain (or semi-infinite due to spectrum symmetry) with data available only on a finite bandwidth. In addition, this procedure does not require direct numerical evaluation of the Hilbert transform. The  Fourier coefficients are  computed by solving an oversampled  regularized least squares problem via a truncated SVD method to have the ill-conditioning of the system under control. Causal Fourier continuations even with moderate number of Fourier coefficients are typically oscillatory in the extended domain but this does not essentially affect the quality of reconstruction of the transfer function on the original frequency domain. The length of the extended domain may be tuned to find more optimal value  that would allow decreasing overall reconstruction errors. The error analysis performed for the proposed method unbiases the error of approximation of a transfer function with a causal Fourier series, the error due to truncation of singular values from the causality violations, i.e. a noise or approximation errors in data obtained from measurements or numerical simulations, respectively. The obtained estimates for upper bounds of these errors can be used to verify causality of the given frequency responses.
The method is applicable to both baseband and bandpass regimes and does not require data points to be equally spaced. It shows high accuracy and robustness 
and is capable of detecting very small localized causality violations of the amplitude close to the machine precision.
The proposed technique is applied to several analytic and simulated examples with and without causality violations. The results demonstrate an excellent performance of the method in agreement with obtained error estimates.

\section{Acknowledgments}

We thank the reviewers for their very detailed and constructive comments. This work was funded by the Micron Foundation. The author L.L.B. would also  like to acknowledge the availability of computational resources made possible through the National Science Foundation Major Research Instrumentation Program: grant no. 1229766.

\bibliography{references_Fourier_continuation}

\begin{thebibliography}{10}
\providecommand{\url}[1]{#1}
\csname url@samestyle\endcsname
\providecommand{\newblock}{\relax}
\providecommand{\bibinfo}[2]{#2}
\providecommand{\BIBentrySTDinterwordspacing}{\spaceskip=0pt\relax}
\providecommand{\BIBentryALTinterwordstretchfactor}{4}
\providecommand{\BIBentryALTinterwordspacing}{\spaceskip=\fontdimen2\font plus
\BIBentryALTinterwordstretchfactor\fontdimen3\font minus
  \fontdimen4\font\relax}
\providecommand{\BIBforeignlanguage}[2]{{%
\expandafter\ifx\csname l@#1\endcsname\relax
\typeout{** WARNING: IEEEtran.bst: No hyphenation pattern has been}%
\typeout{** loaded for the language `#1'. Using the pattern for}%
\typeout{** the default language instead.}%
\else
\language=\csname l@#1\endcsname
\fi
#2}}
\providecommand{\BIBdecl}{\relax}
\BIBdecl

\bibitem{Swaminathan_Engin_2007}
M.~Swaminathan and E.~Engin, \emph{Power Integrity Modeling and Design for
  Semiconductors and Systems}.\hskip 1em plus 0.5em minus 0.4em\relax Prentice
  Hall, 2007.

\bibitem{Gustavsen_Semlyen_1999}
B.~Gustavsen and A.~Semlyen, ``{Rational approximation of frequency domain
  responses by vector fitting},'' \emph{{IEEE Trans. Trans. Power Delivery}},
  vol.~{14}, no.~{3}, pp. {1052--1061}, {1999}.

\bibitem{Deschrijver_Haegeman_Dhaene_2007}
D.~Deschrijver, B.~Haegeman, and T.~Dhaene, ``{Orthonormal vector fitting: A
  robust macromodeling tool for rational approximation of frequency domain
  responses},'' \emph{{IEEE Trans. Adv. Packag.}}, vol.~{30}, no.~{2}, pp.
  {216--225}, {2007}.

\bibitem{Charest_Achar_Nakhla_Erdin_2009}
A.~Charest, R.~Achar, M.~Nakhla, and I.~Erdin, ``{Delay extraction-based
  passive macromodeling techniques for transmission line type interconnects
  characterized by tabulated multiport data},'' \emph{{Analog Integr. Circ.
  S.}}, vol.~{60}, no. {1--2}, pp. {13--25}, {2009}, {50th Midwest Symposium on
  Circuits and Systems, Montreal, CANADA, SEP 05-AUG 08, 2007-2008}.

\bibitem{Triverio_Grivet_Talocia_Nakhla_Canavero_Achar_2007}
P.~Triverio, S.~Grivet-Talocia, M.~S. Nakhla, F.~G. Canavero, and R.~Achar,
  ``{Stability, causality, and passivity in electrical interconnect models},''
  \emph{{IEEE Trans. Adv. Packag.}}, vol.~{30}, no.~{4}, pp. {795--808},
  {2007}.

\bibitem{Oppenheim_Schafer_1989}
A.~V. Oppenheim and R.~W. Schafer, \emph{Discrete-Time Signal Processing}, ser.
  Prentice-Hall Signal Processing Series.\hskip 1em plus 0.5em minus
  0.4em\relax Prentice Hall, 1989.

\bibitem{Blais_Cimmino_Ross_Granger_2009}
J.-F. Blais, M.~Cimmino, A.~Ross, and D.~Granger, ``{Suppression of time
  aliasing in the solution of the equations of motion of an impacted beam with
  partial constrained layer damping},'' \emph{{J. Sound Vib.}}, vol. {326}, no.
  {3-5}, pp. {870--882}, {2009}.

\bibitem{Granger_Ross_2009}
D.~Granger and A.~Ross, ``{Effects of partial constrained viscoelastic layer
  damping parameters on the initial transient response of impacted cantilever
  beams: Experimental and numerical results},'' \emph{{J. Sound Vib.}}, vol.
  {321}, no. {1--2}, pp. {45--64}, {2009}.

\bibitem{Gottlieb_Shu_1997}
D.~Gottlieb and C.-W. Shu, ``{On the {G}ibbs phenomenon and its resolution},''
  \emph{{SIAM Rev.}}, vol.~{39}, no.~{4}, pp. {644--668}, 1997.

\bibitem{Gelb_Tanner_2006}
A.~Gelb and J.~Tanner, ``{Robust reprojection methods for the resolution of the
  {G}ibbs phenomenon},'' \emph{{Appl. Comput. Harmon. Anal.}}, vol.~{20},
  no.~{1}, pp. {3--25}, 2006.

\bibitem{Tadmor_2007}
E.~Tadmor, ``{Filters, mollifiers and the computation of the {G}ibbs
  phenomenon},'' \emph{{Acta Numer.}}, vol.~{16}, pp. {305--378}, 2007.

\bibitem{Mhaskar_Prestin_2009}
H.~Mhaskar and J.~Prestin, ``{Polynomial operators for spectral approximation
  of piecewise analytic functions},'' \emph{{Appl. Comput. Harmon. Anal.}},
  vol.~{26}, pp. {121--142}, 2009.

\bibitem{Beylkin_Monzon_2009}
G.~Beylkin and L.~Monz\'on, ``{Nonlinear inversion of a band-limited {F}ourier
  transform},'' \emph{{Appl. Comput. Harmon. Anal.}}, vol.~{27}, pp.
  {351--366}, 2009.

\bibitem{Kramers_1927}
H.~A. Kramers, ``{La diffusion de la lumiere par les atomes},'' \emph{{Atti
  Cong. Intern. Fisica (Transactions of Volta Centenary Congress) Como}},
  vol.~{2}, pp. {545--557}, {1927}.

\bibitem{Kronig_1926}
R.~D.~L. Kronig, ``{On the theory of dispersion of x-rays},'' \emph{{J. Opt.
  Soc. Am.}}, vol.~{12}, no.~{6}, pp. {547--557}, {1926}.

\bibitem{Guillemin_1977}
E.~A. Guillemin, \emph{Synthesis of Passive Networks: Theory and Methods
  Appropriate to the Realization and Approximation Problems}.\hskip 1em plus
  0.5em minus 0.4em\relax New York: R. E. Krieger, 1977.

\bibitem{Amari_Gimersky_Bornemann_1995}
S.~Amari, M.~Gimersky, and J.~Bornemann, ``{Imaginary part of antennas
  admittance from its real part using {B}odes integrals},'' \emph{{IEEE Trans.
  Antennas Propag.}}, vol.~{43}, no.~{2}, pp. {220--223}, {1995}.

\bibitem{Tesche_1992}
F.~M. Tesche, ``{On the use of the {H}ilbert transform for processing measured
  {C}{W} data},'' \emph{{IEEE Trans. Electromagn. Compat.}}, vol.~{34}, no. {3,
  1}, pp. {259--266}, {1992}.

\bibitem{Knockaert_Dhaene_2008}
L.~Knockaert and T.~Dhaene, ``{Causality determination and time delay
  extraction by means of the eigenfunctions of the Hilbert transform},'' in
  \emph{{2008 IEEE Workshop on Signal Propagation on Interconnects}}, {2008},
  pp. {19--22}, {12th IEEE Workshop on Signal Propagation on Interconnects,
  Avignon, France, May 12-15, 2008}.

\bibitem{Narayana_Rao_Adve_Sarker_Vannicola_Wicks_Scott_1996}
S.~Narayana, G.~Rao, R.~Adve, T.~Sarker, V.~Vannicola, M.~Wicks, and S.~Scott,
  ``{Interpolation/extrapolation of frequency domain responses using the
  Hilbert transform},'' \emph{{IEEE Trans. Microw. Theory Techn.}}, vol.~{44},
  no. {10, 1}, pp. {1621--1627}, {1996}.

\bibitem{Luo_Chen_2005}
S.~P. Luo and Z.~Z. Chen, ``{Iterative methods for extracting causal
  time-domain parameters},'' \emph{{IEEE Trans. Microw. Theory Techn.}},
  vol.~{53}, no. {3, 1}, pp. {969--976}, {2005}.

\bibitem{Young_2010}
B.~Young, ``{Bandwidth and density reduction of tabulated data using causality
  checking},'' in \emph{{2010 IEEE Electrical Design of Advanced Packaging and
  Systems Symposium (EDAPS 2010)}}, {2010}, pp. {1--4}.

\bibitem{Triverio_Grivet_Talocia_2006}
P.~Triverio and S.~Grivet-Talocia, ``{A robust causality verification tool for
  tabulated frequency data},'' in \emph{{10th IEEE Workshop On Signal
  Propagation On Interconnects, Proceedings}}, {2006}, pp. {65--68}, {10th IEEE
  Workshop on Signal Propagation on Interconnects, Berlin, Germany, May 09--12,
  2006}.

\bibitem{Triverio_Grivet_Talocia_2006_2}
------, ``{On checking causality of bandlimited sampled frequency responses},''
  in \emph{{PRIME 2006: 2nd Conference on Ph.D. Research In Microelectronic and
  Electronics, Proceedings}}, {Malcovati, P. and Baschirotto, A.}, Ed., {2006},
  pp. {501--504}, {2nd Conference on Ph.D. Research in MicroElectronics and
  Electronics, Otranto, Italy, June 12-15, 2006}.

\bibitem{Triverio_Grivet_Talocia_2008}
------, ``{Robust causality characterization via generalized dispersion
  relations},'' \emph{{IEEE Trans. Adv. Packag.}}, vol.~{31}, no.~{3}, pp.
  {579--593}, {2008}.

\bibitem{Mandrekar_Swaminathan_2005_2}
R.~Mandrekar and M.~Swaminathan, ``{Causality enforcement in transient
  simulation of passive networks through delay extraction},'' in \emph{{Signal
  Propagation on Interconnects, Proceedings}}, {2005}, pp. {25--28}, {9th IEEE
  Workshop on Signal Propagation on Interconnects, Garmisch Partenkirchen,
  Germany, May 10-13, 2005}.

\bibitem{Mandrekar_Swaminathan_2005}
------, ``{Delay extraction from frequency domain data for causal
  macro-modeling of passive networks},'' in \emph{{2005 IEEE International
  Symposium On Circuits And Systems (ISCAS), Vols 1-6, Conference
  Proceedings}}, ser. {IEEE International Symposium on Circuits and Systems},
  {2005}, pp. {5758--5761}, {IEEE International Symposium on Circuits and
  Systems (ISCAS), Kobe, Japan, May 23-26, 2005}.

\bibitem{Mandrekar_Srinivasan_Engin_Swaminathan_2006}
R.~Mandrekar, K.~Srinivasan, E.~Engin, and M.~Swaminathan, ``{Causal transient
  simulation of passive networks with fast convolution},'' in \emph{{10th IEEE
  Workshop on Signal Propagation on Interconnects, Proceedings}}, {2006}, pp.
  {61--64}, {10th IEEE Workshop on Signal Propagation on Interconnects, Berlin,
  Germany, May 09-12, 2006}.

\bibitem{Lalgudi_Srinivasan_Casinovi_Mandrekar_Engin_Swaminathan_Kretchmer_200%
6}
S.~N. Lalgudi, K.~Srinivasan, G.~Casinovi, R.~Mandrekar, E.~Engin,
  M.~Swaminathan, and Y.~Kretchmer, ``{Causal transient simulation of systems
  characterized by frequency-domain data in a modified nodal analysis
  framework},'' in \emph{{Electrical Performance of Electronic Packaging}},
  {2006}, pp. {123--126}, {15th IEEE Topical Meeting on Electrical Performance
  of Electronic Packaging, Scottsdale, AZ, Oct 23-25, 2006}.

\bibitem{Xu_Zeng_He_Han_2006}
B.~S. Xu, X.~Y. Zeng, J.~He, and D.-H. Han, ``{Checking causality of
  interconnects through minimum-phase and all-pass decomposition},'' in
  \emph{{2006 Conference on High Density Microsystem Design and Packaging and
  Component Failure Analysis (HDP `06), Proceedings}}, {2006}, pp. {271--273}.

\bibitem{Dienstfrey_Greengard_2001}
A.~Dienstfrey and L.~Greengard, ``{Analytic continuation, singular-value
  expansions, and {K}ramers-{K}ronig analysis},'' \emph{{Inverse Problems}},
  vol.~{17}, no.~{5}, pp. {1307--1320}, {2001}.

\bibitem{Aboutaleb_Barannyk_Elshabini_Barlow_WMED13}
H.~A. Aboutaleb, L.~L. Barannyk, A.~Elshabini, and F.~Barlow, ``{Causality
  enforcement of {D}{R}{A}{M} package models using discrete Hilbert
  transforms},'' in \emph{{2013 IEEE Workshop on Microelectronics and Electron
  Devices, WMED 2013}}, {2013}, pp. {21--24}.

\bibitem{Barannyk_Aboutaleb_Elshabini_Barlow_IMAPS}
L.~L. Barannyk, H.~A. Aboutaleb, A.~Elshabini, and F.~Barlow, ``{Causality
  verification using polynomial periodic continuations},'' \emph{{J.
  Microelectron. Electron. Packag.}}, (accept.).

\bibitem{Barannyk_Aboutaleb_Elshabini_Barlow_IMAPS2014}
------, ``{Causality Enforcement of High-Speed Interconnects via Periodic
  Continuations},'' in \emph{{The 47th International Symposium on
  Microelectronics, IMAPS 2014, October 14-16, 2014}}, {2014}.

\bibitem{Cooley_Tukey_1965}
J.~W. Cooley and J.~W. Tukey, ``{An algorithm for machine calculation of
  complex {F}ourier series},'' \emph{Math. Comp.}, vol.~{19}, no.~{90}, pp.
  {297--\&}, {1965}.

\bibitem{Nussenzveig_1972}
H.~M. Nussenzveig, \emph{Causality and Dispersion Relations}.\hskip 1em plus
  0.5em minus 0.4em\relax Academic Press, 1972.

\bibitem{Dym_McKean_1985}
H.~Dym and H.~P. McKean, \emph{Fourier Series and Integrals}, ser. Probability
  and Mathematical Statistics.\hskip 1em plus 0.5em minus 0.4em\relax Academic
  Press, 1985.

\bibitem{Beltrami_Wohlers_1966}
E.~J. Beltrami and M.~R. Wohlers, \emph{Distributions and the Boundary Values
  of Analytic Functions}.\hskip 1em plus 0.5em minus 0.4em\relax Academic
  Press, 1966.

\bibitem{Boyd_2002}
J.~P. Boyd, ``{A comparison of numerical algorithms for {F}ourier extension of
  the first, second, and third kinds},'' \emph{{J. Comput. Phys.}}, vol. {178},
  no.~{1}, pp. {118--160}, {2002}.

\bibitem{Bruno_2003}
O.~Bruno, ``{Fast, high-order, high-frequency integral methods for
  computational acoustics and electromagnetics},'' in \emph{{Topics In
  Computational Wave Propagation: Direct and Inverse Problems}}, ser. {Lecture
  Notes In Computational Science and Engineering}, {Ainsworth, M. and Davies,
  P. and Duncan, D. and Martin, P. and Rynne, B.}, Ed., vol.~{31}, {2003}, pp.
  {43--82}.

\bibitem{Bruno_Han_Pohlman_2007}
O.~P. Bruno, Y.~Han, and M.~M. Pohlman, ``{Accurate, high-order representation
  of complex three-dimensional surfaces via Fourier continuation analysis},''
  \emph{{J. Comput. Phys.}}, vol. {227}, no.~{2}, pp. {1094--1125}, {2007}.

\bibitem{Lyon_2012}
M.~Lyon, ``{Sobolev smoothing of {S}{V}{D}-based {F}ourier continuations},''
  \emph{{Appl. Math. Let.}}, vol.~{25}, no.~{12}, pp. {2227--2231}, {2012}.

\bibitem{Lyon_2012a}
------, ``{Approximation error in regularized {S}{V}{D}-based {F}ourier
  continuations},'' \emph{{Appl. Numer. Math.}}, vol.~{62}, no.~{12}, pp.
  {1790--1803}, {2012}.

\bibitem{Weideman_1995}
J.~A.~C. Weideman, ``{Computing the {H}ilbert transform on the real line},''
  \emph{Math. Comp.}, vol.~{64}, no. {210}, pp. {745--762}, {1995}.

\bibitem{Huybrechs_2010}
D.~Huybrechs, ``{On the {F}ourier extension of nonperiodic functions},''
  \emph{{SIAM J. Numer. Anal.}}, vol.~{47}, no.~{6}, pp. {4326--4355}, {2010}.

\bibitem{LAPACK_Users_Guide}
E.~Anderson, Z.~Bai, C.~Bischof, S.~Blackford, J.~Demmel, J.~Dongarra,
  J.~Du~Croz, A.~Greenbaum, S.~Hammarling, A.~McKenney, and S.~D., \emph{LAPACK
  Users' Guide}, ser. 3rd Ed.\hskip 1em plus 0.5em minus 0.4em\relax SIAM,
  1999.

\bibitem{Trefethen_Bau_1997}
L.~N. Trefethen and D.~Bau~III, \emph{Numerical Linear Algebra}.\hskip 1em plus
  0.5em minus 0.4em\relax SIAM: Society for Industrial and Applied Mathematics,
  1997.

\bibitem{Cheney_2000}
E.~W. Cheney, \emph{Introduction to Approximation Theory}.\hskip 1em plus 0.5em
  minus 0.4em\relax AMS Chelsea Publishing, 2000.

\end{thebibliography}

\begin{IEEEbiography}
[{\includegraphics[width=1in,height=1.25in,clip,keepaspectratio]{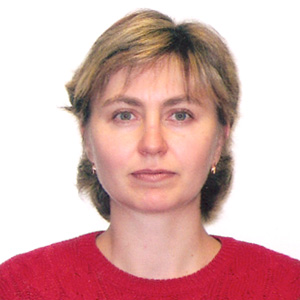}}]{Lyudmyla L. Barannyk}
received the  M.S. degree in Applied Mathematics and Ph.D. degree in Mathematical Sciences from 
New Jersey Institute of Technology in 2000 and 2003, respectively. 

She has been an Assistant Professor in the Department of Mathematics, University of Idaho, since 2007. She was a Postdoctoral Assistant Professor in the Department of Mathematics, University of Michigan, from 2003 to 2007.

Her research interests are electrical modeling and characterization of interconnect packages, scientific computing, mathematical modeling, dimension reduction of large ODE systems, numerical methods for ill-posed problems, fluid dynamics, interfacial instability, PDEs, pseudo-spectral methods, boundary integral methods, grid-free numerical methods.
\end{IEEEbiography}
%
%
\begin{IEEEbiography}
[{\includegraphics[width=1in,height=1.25in,clip,keepaspectratio]{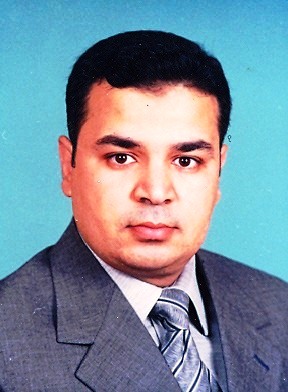}}]{Hazem A. Aboutaleb}
is a Ph.D. student in the Department of Electrical and Computer Engineering, University of Idaho. He received his B.Sc. and M.S. degrees in electrical engineering from The Military Technical College, Cairo, Egypt, in 1994 and 2007, respectively. He was awarded a doctoral fellowship by the Egyptian Government, in 2011, to join University of Idaho as a doctoral student. His research interests are in the area of Microelectronics Fabrication with emphasis on macromodeling of microelectronics package, causality verification and enforcement of macromodels and PI and SI co-simulation analysis.
\end{IEEEbiography}
%
%
\begin{IEEEbiography}
[{\includegraphics[width=1in,height=1.25in,clip,keepaspectratio]{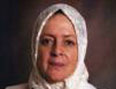}}]
{Aicha Elshabini}
is a Distinguished Professor, Electrical and Computer Engineering department at University of Idaho.  She is a Professional Engineer and IEEE Fellow (CPMT) and IMAPS Fellow.  She is the current advisor for the National Society of Black Engineers (NSBE), the Society of  Women Engineers (SWE), and the International Society of Microelectronics and Electronic Packaging (IMAPS).  Prior to this role, she served as the Dean of Engineering at University of Idaho and the department Head for Electrical Engineering department, and Computer Science and Computer Engineering department at University of Arkansas. 

Professor Elshabini has a B.S.E.E. degree, Cairo University in Electronics \& Communications, a five years academic program (British system), an M.S.E.E., University of Toledo, Ohio in Microelectronics, and a Ph.D., Electrical Engineering from the University of Colorado at Boulder in Solid State Devices and Optoelectronics, 1978.

Elshabini was awarded the 1996 John A. Wagnon Technical Achievements Award, the 2006 Daniel C. Hughes, Jr. Memorial Award,  the 2007 Outstanding Educator Award, and the  2011 President Award.
\end{IEEEbiography}
\begin{IEEEbiography}
[{\includegraphics[width=1in,height=1.25in,clip,keepaspectratio]{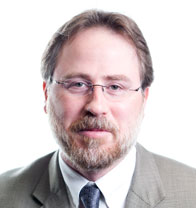}}]{Fred Barlow}
is  a professor and the department chair in the Department  of Electrical and Computer Engineering at the University of Idaho with an emphasis on electronic packaging.  In the past, Dr. Barlow worked for several universities, including Virginia Tech and the University of Arkansas where he held the position of associate department head.

Professor Barlow has served as the major professor for more than twenty graduate students and served as PI or Co-PI for over six million dollars of funded research engaging the microelectronics industry.  Dr. Barlow has published over 100 papers and is coeditor of The Handbook of Thin FilmTechnology (McGraw Hill, 1998), as well as for the Handbook of Ceramic Interconnect Technology (CRC Press, 2007).  He is also a fellow member of the International Microelectronics and Packaging Society (IMAPS) and a senior member of the Institute of Electrical and Electronic Engineers (IEEE). 
\end{IEEEbiography}

\vfill

\end{document}